\documentclass[12pt]{amsart}
\usepackage{amsmath,amstext,amssymb,amscd}
\usepackage{pb-diagram}
\usepackage{fullpage}
\usepackage{eucal}
\usepackage{verbatim}

\theoremstyle{plain}
\swapnumbers
\newtheorem{thm}{Theorem}[section]
\newtheorem{cor}[thm]{Corollary}
\newtheorem{prop}[thm]{Proposition}

\theoremstyle{definition}
\newtheorem{example}[thm]{Example}
\newtheorem{defn}[thm]{Definition}
\newtheorem{thirdarrowproperty}[thm]{Third Arrow Property}
\newtheorem*{question}{Question}
\newtheorem*{remark}{Remark}

     \newcommand{\sB}{\mathcal B}
     \newcommand{\sC}{\mathcal C}
     \newcommand{\sD}{\mathcal D}
     \newcommand{\sE}{\mathcal E}
     
     \newcommand{\sG}{\mathcal G}
     \newcommand{\sH}{\mathcal H}
     
     \newcommand{\sK}{\mathcal K}

     \newcommand{\sO}{\mathcal O}
     \newcommand{\sP}{\mathcal P}

\def\bbC{\mathbb C}

\def\bbN{\mathbb N}

\def\bbR{\mathbb R}
\def\bbT{\mathbb T}
\def\bbZ{\mathbb Z}

\newcommand{\Id}{\mbox{\rm Id}}
\newcommand{\cstar}{\ensuremath{\text{C}^{*}}-}
\renewcommand{\star}{\ensuremath{{}^{*}}\nobreakdash-\hspace{0 pt}}
\DeclareMathOperator{\dom}{dom}
\DeclareMathOperator{\ran}{ran}
\DeclareMathOperator{\supp}{supp}

\begin{document}

\title{Analytic Partial Crossed Products}

\author[A.P. Donsig]{Allan P. Donsig}

\address{Dept. of Mathematics \& Statistics\\
         University of Nebraska---Lincoln\\
         Lincoln, NE,  68588}
\email{adonsig@math.unl.edu}
\author[A. Hopenwasser]{Alan Hopenwasser}
\address{Dept. of Mathematics\\
        University of Alabama\\
        Tuscaloosa, AL 35487}
\email{ahopenwa@euler.math.ua.edu}
\date{May 19, 2003}
\begin{abstract}
Partial actions of discrete abelian groups can be
used to construct
both groupoid \ensuremath{\text{C}^{*}}-algebras and
partial crossed product algebras.  In each case there is a natural
notion of an analytic subalgebra.  We show that for countable
subgroups of $\mathbb R$ and free partial actions, these
constructions yield the same
 \ensuremath{\text{C}^{*}}-algebras and the same analytic
subalgebras.

We also show that under suitable hypotheses an analytic
partial crossed product preserves all the information
in the dynamical system in the sense that two analytic
partial crossed products are isomorphic as Banach algebras
if, and only if, the partial actions are conjugate.
\end{abstract}

\maketitle

\section{Introduction} \label{s:intro}

One way to obtain interesting non-self-adjoint subalgebras of
\cstar algebras is to restrict a crossed product to the subalgebra
associated with a positive cone in the group.  A significant
limitation to this approach arises from the fact that very important
\cstar algebra contexts, for example, AF \cstar algebras, do not
appear as crossed products.  In \cite{MR95g:46122}, Exel found a new
construction, the partial crossed product, that gives many of these 
contexts, by using partial actions.
Exel considered partial actions by $\bbZ$, the integers; 
in~\cite{MR96i:46083} McClanahan extended these ideas to
discrete groups.

Finite dimensional \cstar algebras and AF \cstar algebras can be
realized as partial actions of abelian ordered groups on an abelian
 \cstar algebra (or, equivalently, by a partial action on the
spectrum of the abelian \cstar algebra).  This allows us to specify
``analytic'' subalgebras of the partial crossed product by restricting
to the subalgebra generated by the positive cone in the group.
To indicate the range of algebras that may be obtained by this
construction, we realize Power's toroidal limit 
algebras~\cite{MR93d:46092}
as analytic partial crossed products (Example~\ref{toroidal}).
We should point out that the \cstar envelopes of these algebras are
the Bunce-Deddens \cstar algebras.

AF \cstar algebras are all groupoid \cstar algebras; as such they also
possess analytic subalgebras which are defined in terms of cocycles on
the groupoid.  This raises the question: how do analytic partial
crossed products and analytic subalgebras of groupoid \cstar algebras
fit together?  One of the main results in this paper,
Theorem~\ref{T:bicor} addresses this
issue: free partial actions of countable
discrete subgroups of $\bbR$ on separable abelian \cstar algebras
yield the same \cstar algebras as do locally compact $r$-discrete
principal groupoids with second countable unit space and a locally
constant cocycle.  Most importantly, the analytic subalgebras in the
two constructions coincide.  From a slightly different perspective,
partial actions lead to two constructs: a partial crossed product and a
groupoid \cstar algebra, each with a natural analytic subalgebra.
Theorem \ref{T:bicor} says that these two constructs, together with
the analytic subalgebras, coincide.  This makes available the theory
of partial crossed products for the study of many non-self-adjoint
algebras which have been primarily investigated via groupoid
techniques.

One of the motivations for the study of analytic crossed products is
that, unlike \cstar algebra crossed products, they preserve the
information contained in the dynamical systems which determine the
crossed product.  Results of this type can be found
in~\cite{MR35:1751,MR40:3322,MR86e:46063,MR89e:47069,%
MR90a:46175,MR93h:46095}.

In Section \ref{s:conjugacy} we extend to partial actions of
$\bbZ^+$ a theorem of Power \cite{MR93h:46095} for analytic crossed
products.  Our result states that if $\alpha$ and $\beta$ are free
partial actions of $\bbZ^+$ with the property that the domains of
$\alpha_n$ and $\beta_n$ shrink as $n$ increases in $\bbZ^+$, then
$C_0(X) \times_{\alpha} \bbZ^+$ and
$C_0(X) \times_{\beta} \bbZ^+$ are isomorphic as Banach algebras if,
and only if, the partial actions $\alpha$ and $\beta$ are conjugate in
a natural sense.

While we are primarily interested in analytic crossed products of the
form $C_0(X) \times_{\alpha} \Sigma$, where $\Sigma$ is a positive
cone in a group $G$, these are defined in terms of the \cstar algebra
partial crossed product $C_0(X) \times_{\alpha} G$.
As a first step in the direction of developing a theory which does not
depend on the whole group and the \cstar algebra, we modify the
definition of partial action to a form suitable for a positive cone
in a group.
If our axioms for a partial action of $\Sigma$ are applied to a group,
they simply form a redundant version of McClanahan's definition.
We show by example that these conditions are not redundant when
applied to a cone $\Sigma$.

The question of the extension of a partial action from $\Sigma$ to the
group $G$ arises immediately.
This is solved easily when $\Sigma$ totally orders $G$.
On the other hand, we give an example of a partial action of $(\bbZ^3)^+$
on a seven point space (or, equivalently, on a seven dimensional abelian
\cstar algebra)
which has no extension to $\bbZ^3$.
This material is presented in Section \ref{s:extend}.

Theorem~\ref{T:bicor} in Section \ref{s:anco} makes substantial use of
groupoids.  For treatises with extensive treatment of groupoids, we
direct the reader to \cite{MR82h:46075,MR2001a:22003}.
A helpful brief introduction can be found in \cite{MR90m:46098}.
For the convenience of the reader, we record a few definitions here.
A groupoid is a set $\sG$ with a partially defined multiplication and
an inversion.
If $a$ and $b$ can be multiplied, then $(a,b)$ is called a
\emph{composable pair}; $\sG^2$ denotes the set of all composable pairs.
The multiplication satisfies an associative law and inversion satisfies
$(a^{-1})^{-1}=a$.
Elements of the form $a^{-1}a$ and $aa^{-1}$ are called \emph{units}.
Units act as left and right identities when multiplied by elements with
which they are composable.

We shall only be interested in \emph{principal groupoids}; these are
all equivalence relations.  If $\sG$ is an equivalence relation on a
set $X$, then $(x,y)$ and $(w,z)$ are composable if, and only if,
$y=w$.  In this case, $(x,y)(y,z)=(x,z)$.  Also
$(x,y)^{-1}=(y,x)$, always.  Groupoids have \emph{range} and
\emph{domain} maps defined by $r(a) = aa^{-1}$ and
$d(a) = a^{-1}a$.  In the principal groupoid context, this gives
$r(x,y)=(x,x)$ and $d(x,y)=(y,y)$.  Since we can identify the diagonal
$\{(x,x) \mid x \in X\}$ with $X$ in a natural way, we can view $r$ and
$d$ as maps from $\sG$ onto $X$; they are, in fact, just the
coordinate projections.  A subset $E \subseteq \sG$ is said to be a
$\sG$-\emph{set} if $r$ and $d$ are both one-to-one on $E$.  In this
case, $E$ is just the graph of a function from a subset of $X$ to some
other subset.

All the groupoids which we consider will carry a
 locally compact topology;
the groupoid operations will, of course, be continuous with respect to
this topology.
An $r$-\emph{discrete} groupoid is one in which the 
set of units is open.
In our context (principal groupoids) this means that the diagonal
$\{(x,x) \mid x \in X\}$ is an open subset of $\sG$.
It will, in fact, be the case that the open $\sG$-sets form a base
for the topology.
One consequence of note is that in $r$-discrete principal groupoids,
each equivalence class is countable.  Finally, a \emph{cocycle} on
$\sG$ is a continuous map $c \colon \sG \to \bbR$ such that
$c(x,z) = c(x,y) + c(y,z)$, for all composable pairs $(x,y)$ and
$(y,z)$ in $\sG$.

\section{Partial Actions} \label{s:partial}

The following definition is an extension, suitable for use in the
context of non-selfadjoint algebras, of the definition of a
partial action of an abelian group on a \cstar algebra.
For the usual definition, see \cite{MR95g:46122,MR96i:46083}.

\begin{defn}
Let $A$ be a \cstar algebra, let $(G,+)$ be
 a discrete abelian group,
and let $\Sigma$ be a subset of $G$ for which
$\Sigma + \Sigma \subset \Sigma$.
A \emph{partial action} of $\Sigma$ on $A$ is a family of isomorphisms
$\alpha = \{ \alpha_t \colon D_{-t} \to D_t \mid t \in \Sigma \}$
between closed two-sided ideals of $A$ so that

\begin{description}
\item[1a] $\alpha_s(D_{-s} \cap D_t) = D_s \cap  D_{s+t}$,
  \quad if $s,t \in \Sigma$,
\item[1b] $\alpha_t(D_{-t} \cap D_{-s-t}) = D_t \cap  D_{-s}$,
  \quad if $s,t \in \Sigma$,
\item[2] $\alpha_{s+t}(x) = \alpha_s \circ \alpha_t(x)$,
  \quad if $x \in D_{-t} \cap D_{-s-t}$,\quad  and
\item[3] $D_0=A$ and $\alpha_0 = \Id_A$.
\end{description}
\end{defn}

We will show below that Conditions 1a, 1b, and 2 are equivalent to
a property which is usually much easier to work with.
We shall refer to this property as the \emph{third arrow property}.

\begin{thirdarrowproperty}
Any two of the following statements
implies the third.
\begin{description}
\item[I] $x \in D_{-t}$ and $y = \alpha_t(x)$.
\item[II] $y \in D_{-s}$ and $z  = \alpha_s(y)$.
\item[III] $x \in D_{-s-t}$ and $z = \alpha_{s+t}(x)$.
\end{description}
\end{thirdarrowproperty}

This can be indicated schematically by saying that if two of the arrows
in the following diagram exist,
then so does the third arrow.

\[
\begin{diagram}
\node{x}\arrow[1]{ese,b} {\alpha_t}   \arrow[4]{e,t}{\alpha_{s+t}}
 \node[4]{z} \\
\node[3]{y} \arrow[1]{ene,b}{\alpha_s}
\end{diagram}
\]

Another way of saying this is that in the order on the orbit of a point
induced by the partial action, any two points are comparable.
(The order is a total order.)

One consequence of the third arrow property (of I and II implies III,
to be specific) is that for all $s,t \in \Sigma$,
$\alpha_s \circ \alpha_t$ is a restriction of $\alpha_{s+t}$.

Moreover, when $\Sigma$ is the whole group,
 either of  conditions 1a or 1b is also equivalent to
another condition, which is weaker in the general case.

\begin{description}
\item[1$^\prime$] $\alpha_s(D_{-s} \cap D_t) \subset D_{s+t}$,
\quad if $s,t \in \Sigma$,
\end{description}

Condition $1'$ says exactly that statements I and II in the third
arrow property imply
statement III.  In the case when $\Sigma$ is a group,
all arrows can be reversed; so each of the three
implications in the third arrow property implies the other two.
Conditions $1'$, 2 and 3 are exactly the definition of partial action
of a group, as given by McClanahan in \cite{MR96i:46083}.

\begin{proof}[Proof of Equivalence of the Third Arrow Property.]
First, assume that the three conditions in the definition hold.
We only show that Conditions II and III imply I;
similar arguments apply to the other implications.
Schematically
\[
\begin{diagram}
\node{x}\arrow[1]{ese,b,..} {\alpha_t}   \arrow[4]{e,t}{\alpha_{s+t}}
 \node[4]{z} \\
\node[3]{y} \arrow[1]{ene,b}{\alpha_s}
\end{diagram}
\]
Since $z \in D_s \cap D_{s+t}$, condition 1a implies that there exists
$y' \in D_{-s} \cap D_t$ such that $\alpha_s(y') = z$.  But
$y \in D_{-s}$, $\alpha_s(y) = z$, and $\alpha_{-s}$ is injective on
$D_{-s}$; this yields $y' = y$.  Now we know that
 $y \in D_t \cap D_{-s}$.  Condition 1b gives the existence of
$x' \in D_{-t} \cap D_{-s-t}$ such that $\alpha_t(x') = y$.  By
condition 2,
$z=\alpha_s(y)=\alpha_s(\alpha_t(x'))=\alpha_{s+t}(x')$.
Thus, $x,x' \in D_{-s-t}$ and
$\alpha_{s+t}(x) = \alpha_{s+t}(x')$.  Since $\alpha_{s+t}$ is
injective on $D_{-s-t}$, $x=x'$.  This yields
$\alpha_t(x)=y$ and condition II.

For the converse, assume that the third arrow property holds.
We first verify Condition~1b; the argument for Condition~1a is similar.
Assume that $x \in D_{-t} \cap D_{-s-t}$.
Then there exist $y \in D_t$ such that $\alpha_t(x) =y$ and
$z \in D_{s+t}$ such that $\alpha_{s+t}(x)=y$.  This is I and III, so
II holds and $y \in D_{-s}$.  This gives
$\alpha_t(D_{-t} \cap D_{-s-t}) \subseteq D_t \cap D_{-s}$.

Now let $y \in D_t \cap D_{-s}$.  Then there exist $x \in D_{-t}$ such
that $\alpha_t(x) = y$ and $z \in D_s$ such that
$\alpha_s(y) = z$.  This is I and II; III now gives
$x \in D_{-s-t}$ and
$y = \alpha_t(x) \in \alpha_t(D_{-t} \cap D_{-s-t})$.  Thus
$D_t \cap D_{-s} \subset \alpha_t(D_{-t} \cap D_{-s-t})$ and
condition 1b is verified.

Condition 2 follows in a similar way from the statement that
I and III imply II.
\end{proof}

Most of the time, we shall assume that $\Sigma \subset G$ is a cone,
i.e., $\Sigma + \Sigma \subset \Sigma$
and $\Sigma \cap -\Sigma = \{0\}$.
We will be most interested in
the context $A=C_0(X)$, where $X$ is a locally
compact Hausdorff space and $G$ is a discrete subgroup of $(\bbR,+)$
with $\Sigma=G \cap [0,+\infty)$.

Partial actions on a locally compact Hausdorff space $X$ are defined
in the same way as partial actions on \cstar algebras, except that the
ideals are replaced by open subsets of $X$ and the \star isomorphisms
between ideals are replaced by homeomorphisms between open sets.
Often, it is convenient to ``move'' a partial action from $C_0(X)$ to
$X$.  This is routine, but here is a sketch of the procedure.

If $D_t$ is an ideal in $C_0(X)$, then there is a closed subset $C_t$
of $X$ such that
$D_t = \{f \in C_0(X) \mid f \equiv 0 \text{ on } C_t \}$.  Then
$X_t = X \setminus C_t$ is an open subset of $X$ and
$D_t \cong C_0(X_t)$.  With this notation, given an isomorphism
$\alpha_t \colon C_0(X_{-t}) \to C_0(X_t)$, there is a
homeomorphism $\beta_t \colon X_{-t} \to X_t$ such
that, for all $f \in  C_0(X_{-t})$,
$\alpha_t(f) = f \circ \beta_t^{-1}$.
It is straightforward to check that
$\alpha$ is a partial action on $C_0(X)$ if, and only if, $\beta$ is a
partial action on $X$.

The first two examples below are closely related to the standard 
and refinement triangular subalgebras of UHF \cstar algebras.

\begin{example} \label{eg:std}
Let $X$ be the Cantor set $\prod_{i=1}^\infty \{0,1\}$ and let 
$\alpha_1$ be the partial map on $C(X)$ induced by the odometer 
map $\beta$: if $x \in X$ is not $(1,1,\ldots)$, then $\beta(x)$ 
is given by adding $1$ to the first coordinate of $x$, with carries 
to the right.
Then $\alpha_n$ is defined to be $(\alpha_1)^n$ on the domain where 
this makes sense.  
This $\alpha$ defines a partial action of $\Sigma = \bbZ^+$ on $C(X)$.
\end{example}

\begin{example} \label{eg:refine}

Let $G$ be the dyadic rationals
$\{ k/2^n \mid k \in \bbZ, n \in \bbN \}$ and
let $\Sigma = G \cap [0,+\infty)$.
Let $X =\prod_{i=1}^\infty \{0,1\}$, the Cantor set, where we
associate elements of $X$ with base $2$ representations of numbers
in $[0,1]$.
Note that dyadic rationals in $[0,1]$ have two expansions so,
for example, $1/2$
becomes two numbers: $1/2^+=.1000\ldots$ and $1/2^-=.0111\ldots$.
For $s \in \Sigma$ we define $\beta_s(x)$ to be $x+s$, provided
$x+s$ is in $[0,1]$.
Thus, $\beta_{1/2}$ is defined for $\{ (x_i) \in X \mid x_1=0 \}$
and sends $(0,x_2,x_3,\ldots)$ to $(1,x_2,\ldots)$.
Then $\beta$ is a partial action of $\Sigma$ on $X$.
\end{example}

\begin{example} \label{eg:split1a1b}
Neither of conditions 1a and 1b implies the other.  This is
illustrated by modifications of a simpler 
version of Example~\ref{eg:refine}.
The modifications are not partial actions, of course.

Let $X = (0,1)$ and $\Sigma = [0, \infty) \subseteq \bbR$.  For each
$t \in \Sigma$, let
\[
X_{-t} =
\begin{cases}
(0, 1-t), &\text{if } t < 1 \\
\emptyset, &\text{if } t \geq 1
\end{cases}
\]
and
\[
X_t =
\begin{cases}
(t,1), &\text{if } t < 1 \\
\emptyset &\text{if } t \geq 1.
\end{cases}
\]
Define $\beta_t \colon X_{-t} \to X_t$ by $\beta_t(x) = x+t$. 
It is easy to see that $\beta$ satisfies the third arrow property 
and Condition 3.
If countable groups are desired, restrict $\alpha$ to the intersection
of $\Sigma$ with any countable, dense subgroup of $\bbR$.

We now modify $\beta$ to obtain a system of partial homeomorphisms
of $X$ which satisfies condition 1a but not condition 1b.  We take
$\gamma_0 = id_X \;(= \beta_0)$ and, for $t \geq 1$,
$\gamma_t = \emptyset \;(= \beta_t)$.  For $0<t<1$, the values of
interest,
 define $\gamma$ as follows:
\begin{align*}
&\text{for } 0 < t \leq 1/2, &&\gamma_t \text{ is }
\beta_t \text{ restricted to } (0,1/2) = D_{-t} \setminus
[1/2,1-t], \\
&\text{for } 1/2 \leq t <1, &&\gamma_t = \beta_t.
\end{align*}
That is, $\gamma$ is obtained from $\beta$ by taking
$G=\{(x,y) \mid 0<x \leq y <1\}$, the union of the graphs of the
$\beta_t$ and deleting the subset
$\{(x,y) \in G \mid x \geq 1/2 \text{ and } x \neq y\}$.

This example satisfies I and II implies III and also
II and III implies I of the third arrow property; it does not satisfy
I and III implies II.  To be specific, take $x=3/8$, $y=5/8$,
$z=6/8$, $t=2/8$ and $s=1/8$.  Then
\begin{align*}
&\gamma_t(x)=\gamma_{2/8}(3/8)=5/8=y \text{ and} \\
&\gamma_{s+t}(x)=\gamma_{3/8}(3/8)=6/8=z,
\end{align*}
but $y=5/8$ is not in the domain of $\gamma_s=\gamma_{1/8}$, so
we fail to have $\gamma_s(y)=z$.

In terms of the conditions in the definition of a partial action,
$\beta$ satisfies 1a but not 1b.

If instead, we define $\gamma$ on the essential values by
\begin{align*}
&\text{for } 0 < t \leq 1/2, &&\gamma_t \text{ is } \beta_t
\text{ restricted to } (1/2-t,1-t) = D_{-t} \setminus (0,1/2-t],\\
&\text{for } 1/2 \leq t <1, &&\gamma_t = \beta_t,
\end{align*}
then we obtain an example in which I and II implies III and
I and III implies II but II and III do not imply I.  
(Condition 1b is satisfied but Condition 1a is not.)  
This example is obtained by deleting
$\{(x,y) \in G \mid y \leq 1/2 \text{ and } x \neq y\}$ from $G$.

Finally, deletion of
$\{(x,y) \in G \mid y \leq 1/2 \text{ or }
 x \geq 1/2 \text{ and } x \neq y \}$
yields an
example in which the only implication in the third arrow property
to hold is I and II implies III.
Conditions 1a and 1b in the definition both fail.
\end{example}

\begin{defn}
We say that a partial action
is \emph{non-degenerate} if the group
generated by $\{ s \in \Sigma \mid D_s \ne \emptyset \}$ is $G$.
By replacing $G$ with the subgroup generated by this set, we may
always assume a partial action is non-degenerate.
\end{defn}

\begin{defn}
A partial action satisfies the \emph{composition property} if
$\alpha_{s+t} = \alpha_s \circ \alpha_t$
for all  $s,t \in \Sigma$.
The partial action satisfies the
\emph{domain ordering property} if
$D_{-s-t} \subseteq  D_{-t}$, for all $s,t \in \Sigma$.
\end{defn}

\begin{remark}
These two definitions are, in fact, equivalent.
It is trivial to see that the composition property implies
the domain ordering property.
Now assume that the domain ordering condition holds.    We have
already seen that whenever $\alpha$ is a partial action of $\Sigma$,
 $\alpha_s \circ \alpha_t$ is a restriction of
$\alpha_{s+t}$, so it remains
only to show
that $\dom \alpha_{s+t} \subseteq \dom \alpha_s \circ \alpha_t$.
 Let $x \in   D_{-s-t}$.  By the domain ordering property, we
also have $x \in D_{-t}$.  Letting $z = \alpha_{s+t}(x)$ and
$y = \alpha_t(x)$, the third arrow property (I and III implies II)
tells us that $y \in D_{-s}$ and $z = \alpha_s(y)$.  In particular,
$x \in \dom \alpha_s \circ \alpha_t$.
\end{remark}

The next example gives a
partial action which fails to satisfy
these two equivalent conditions.

\begin{example} \label{eg:arc}
  Let $X$ be an arc in the unit
circle; i.e., $X = \{e^{it}\mid a<t<b\}$, for some values of $a$ and
$b$. Define a partial action of
$\bbZ$  on $X$ as follows: for each $n$, let
$X_{-n} =\mbox{$\{x\in X \mid e^{in}x \in X\}$}$ and let
$\alpha_n(x) = e^{in}x$, for  $x\in X_{-n}$.
  Let $\Sigma = \bbZ^+$.  The
restriction of $\alpha$ to $\Sigma$ is a partial action
(the third arrow property is particularly easy to verify)
 which need not satisfy the composition condition.
  In particular, if $a=0$ and $b=1$
 then $X_{-3} = \emptyset$ while
$X_{-7} = \{e^{it} \mid 0 < t < 2\pi-6 \approx .28 \}$.
Hence, $\alpha_4 \circ \alpha_3$ is a proper restriction of
$\alpha_7$.

Note that $X_{-6}=\{e^{it} \mid 2\pi-6 <t<1\} \ne \emptyset$.  Since 6
and 7 are relatively prime, this partial action is non-degenerate.
\end{example}

\begin{example}
Consider $X=\bbR$ and $\Sigma=\bbZ^+$.
For $t$ odd, define $X_t$ to be $\cup_{a \in \bbZ} [a,a+1/2]$
and $\beta_t$ to be translation by $t$.
For $t$ even, define $X_t=\bbR$ and $\beta_t$ to be translation by $t$.
Then the action is non-degenerate but $X_5 \subset X_6$.
We can write this as a `direct sum' of a non-degenerate action on $X_1$ 
and a degenerate action on $\bbR\backslash X_1$.

This example can be generalized by considering $\bbN \times [0,1]$,
$2\bbN \times [0,1]$, $4\bbN \times [0,1]$, and so on.
\end{example}

\section{Extensions} \label{s:extend}

Assume that $G$ is an abelian group and that $\Sigma$ is a subset of
$G$ which satisfies:
\begin{enumerate}
\item $\Sigma + \Sigma \subseteq \Sigma$. \label{pc:closure}
\item $\Sigma \cap -\Sigma = \{0\}$. \label{pc:intersect}
\item $G = \Sigma - \Sigma$.  \label{pc:generate}
\end{enumerate}

A subset which satisfies these properties will be referred to as a
\emph{positive cone} (or, sometimes, simply as a \emph{cone}).
The pair $(G,\Sigma)$ will be called a \emph{directed group}.
 If, in
addition, $G = \Sigma \cup (-\Sigma)$, then we say that $\Sigma$
\emph{totally orders} $G$.  The total ordering is given by
$a \leq b$ if, and only if, $b-a \in \Sigma$.

\begin{question}
Suppose that $\alpha$ is a partial action of a positive cone $\Sigma$
acting on a \cstar algebra $A$.  Is there a (necessarily unique)
extension of $\alpha$ to $G$?
\end{question}

This question is easy to settle in the totally ordered case:

\begin{prop} \label{prop:toextension}
Assume that $G$ is totally ordered by a cone $\Sigma$.
If $\alpha$ is a partial action of $\Sigma$ on a \cstar algebra $A$,
then $\alpha$ has a unique extension to a partial action of $G$ on $A$.
\end{prop}

\begin{proof}
Since $G = \Sigma \cup (-\Sigma)$, we only need to define $\alpha$
appropriately on $-\Sigma$ and it is obvious how this must be done:
for $t \in \Sigma$, take $\alpha_{-t} = \alpha_t^{-1}$.
Condition 3 in the definition of partial action is trivially
satisfied.  Rather than verifying conditions 1a, 1b, and 2, it is
easier to verify the third arrow property.

A ``hand waving'' proof is elementary: given a diagram consisting
of three elements in suitable ideals and two
arrows indexed by elements of $G$, we can, if necessary, reverse one
or both of the arrows to obtain a diagram in which all arrows are
indexed by elements of $\Sigma$.  Since $\alpha$ satisfies the third
arrow property on $\Sigma$, we obtain the third arrow in the new
diagram.  Either this arrow or its inverse will yield the required
third arrow for the original diagram.

Here is a more detailed argument showing that I and III implies II.
The other two implications needed for the third arrow property can be
verified in a similar fashion.
We assume, then, that $x \in D_{-s-t}$,
$x \in D_{-t}$, $z = \alpha_{s+t}(x)$, and $y = \alpha_t(x)$.
\[
\begin{diagram}
\node{x}\arrow[1]{ese,b} {\alpha_t}   \arrow[4]{e,t}{\alpha_{s+t}}
 \node[4]{z} \\
\node[3]{y}
\end{diagram}
\]
We need to deal with the four possibilities regarding the membership
of $s$ and $t$ in $\Sigma$ and $-\Sigma$.

First assume that $t$ and $s+t$ both lie in $\Sigma$.  If $s$ is also
in $\Sigma$, then the third arrow property for $\Sigma$ yields
$y \in D_{-s}$ and $z = \alpha_s(y)$;
 thus II holds.  If, instead, $s \in -\Sigma$, then
view $\alpha_t$ as III and $\alpha_{s+t}$ as I; the third arrow
property for $\Sigma$ yields $z \in D_s$ and
 $y = \alpha_{-s}(z)$.  But then $y \in D_{-s}$ and
$z = \alpha_s(y)$; once again II holds.

Now assume that $t \in -\Sigma$ and $s+t \in \Sigma$.  This forces
$s \in \Sigma$.  Consider the diagram
\[
\begin{diagram}
\node{x}  \arrow[4]{e,t}{\alpha_{s+t}}
 \node[4]{z}  \\
\node[3]{y} \arrow[1]{wnw,b} {\alpha_{-t}}
\end{diagram}
\]
in which both indices are in $\Sigma$.
We can treat
$\alpha_{-t}$ as I and $\alpha_{s+t}$ as II to obtain
$z = \alpha_{-t +s +t}(y) = \alpha_s(y)$.  In the original diagram,
this is II, just what we need to deduce.

Next assume that $t \in \Sigma$ and $s+t \in -\Sigma$.  Then we must
have $s \in -\Sigma$.  The appropriate diagram with indices in
$\Sigma$ is
\[
\begin{diagram}
\node{x}\arrow[1]{ese,b} {\alpha_t}
 \node[4]{z} \arrow[4]{w,t}{\alpha_{-s-t}}\\
\node[3]{y}
\end{diagram}
\]
Use the third arrow property for $\Sigma$ (the I and II implies III
component) to obtain
$y = \alpha_{-s-t+t}(z) = \alpha_{-s}(z)$.  But then
$z = \alpha_s(y)$ and yet again II is verified for the original
diagram.

Finally suppose that $t \in -\Sigma$ and $s+t \in -\Sigma$.  This
gives a diagram
\[
\begin{diagram}
\node{x}
 \node[4]{z} \arrow[4]{w,t}{\alpha_{-s-t}}\\
\node[3]{y} \arrow[1]{wnw,b} {\alpha_{-t}}
\end{diagram}
\]
with indices in $\Sigma$.  By the first part of the detailed proof,
$y = \alpha_{-s}(z)$ and again $z = \alpha_s(y)$.

As mentioned earlier, the other two implications in the third arrow
property can be verified with similar arguments.
\end{proof}

Returning to the general case, to extend $\alpha$ to all of $G$
it would be natural to define $\alpha_g$ to be the
union of all compositions
\[ \alpha_{t_n}^{-1} \circ \alpha_{s_n} \circ \cdots
        \circ \alpha_{t_1}^{-1} \circ \alpha_{s_1} \]
where $s_1,\cdots,s_n,t_1,\cdots,t_n$ are in $\Sigma$
and $g=-t_n+s_n - \cdots -t_1 + s_1$.
The third arrow property for $\Sigma$ implies that we need
only consider expressions for $g$ of this form, as we may
compose adjacent elements with both indices in $\Sigma$ or
both indices in $-\Sigma$.

We next find a sufficient condition for
$\alpha_g$ to be well-defined. Suppose that
$x\in A$ and $g \in G$,  that $g$ has  two expressions
$g=-t_n+s_n \cdots -t_1 + s_1$ and
$g=-v_k+u_k \cdots -v_1 + u_1$, and that $x$ is in the
domain of both the associated compositions.
If $x'= \alpha_{t_n}^{-1} \circ \alpha_{s_n} \circ \cdots
        \circ \alpha_{t_1}^{-1} \circ \alpha_{s_1}(x)$,
then $x'$ is in the domain of the composition associated to
\[
-v_k+u_k \cdots -v_1 + u_1 +t_n-s_n \cdots +t_1 - s_1.
\]
This expression sums to the identity element of the group,
so the associated composition must send $x'$ to $x'$.
If this always holds, then $\alpha_g$ is well-defined; it is easy
to see that the third arrow property now holds.

The condition of the previous paragraph is clearly valid for
a partial action of a group, so this condition is both necessary
and sufficient to have an extension.

Summarizing, we have the following observation.

\begin{prop} \label{p:ext}
Let $(G,\Sigma)$ be a directed group and let
 $\alpha$ be a partial action of $\Sigma$.
Then $\alpha$ extends to a partial action of $G$ if and only if
for all $s_1,\ldots,s_n$ and $t_1,\ldots,t_n$ in $\Sigma$ with
$-t_n+s_n-\ldots-t_1+s_1$ equal to the identity,
\[ \alpha_{t_n}^{-1} \circ \alpha_{s_n} \circ \cdots
        \circ \alpha_{t_1}^{-1} \circ \alpha_{s_1} \]
is a restriction of the identity map.
\end{prop}

The next example shows that there are partial actions which 
do not extend.

\begin{example} \label{no_ext}
Let $X=\{1,2,3,4,5,6,7\}$ with the discrete topology
and $G=\bbZ^3$ with the usual positive cone.
Define $\beta_{(1,0,0)}$ to have domain $\{1,5\}$
and send $1$ to $2$ and $5$ to $4$.
Define $\beta_{(0,1,0)}$ to have domain $\{3,5\}$
and send $3$ to $2$ and $5$ to $6$.
Define $\beta_{(0,0,1)}$ to have domain $\{3,7\}$
and send $3$ to $4$ and $7$ to $6$.

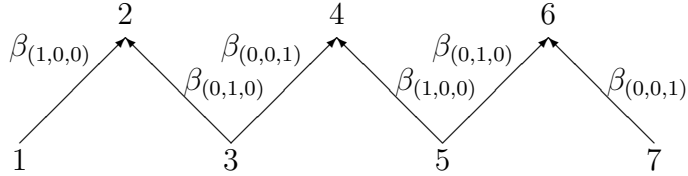
\begin{center}
\begin{picture}(240,120)
\put(  0,40){\vector(1,1){40}}
\put( 80,40){\vector(-1,1){40}}
\put( 80,40){\vector(1,1){40}}
\put(160,40){\vector(-1,1){40}}
\put(160,40){\vector(1,1){40}}
\put(240,40){\vector(-1,1){40}}
\put(-3,30){1}
\put(77,30){3}
\put(157,30){5}
\put(237,30){7}
\put(37,85){2}
\put(117,85){4}
\put(197,85){6}
\put(-4,73){$\beta_{(1,0,0)}$}
\put(62,60){$\beta_{(0,1,0)}$}
\put(76,73){$\beta_{(0,0,1)}$}
\put(142,60){$\beta_{(1,0,0)}$}
\put(156,73){$\beta_{(0,1,0)}$}
\put(222,60){$\beta_{(0,0,1)}$}
\end{picture}
\end{center}

Since the domains of these maps are contained in $\{1,3,5,7\}$
and the ranges are contained in $\{2,4,6\}$,
all possible compositions are empty.
Thus, for all other strictly positive elements of $G$, we may
define $\beta_g$ to be the trivial map with empty domain.
Of course, $\beta_{(0,0,0)}$ is the identity map.
The third arrow property holds trivially, since whenever
there is a point in $X$ that is in the domain of two different
elements of $\Sigma$, the two elements are incomparable.

Since $\beta^{-1}_{(0,0,1)} \circ \beta_{(0,1,0)} \circ
        \beta^{-1}_{(1,0,0)} \circ \beta_{(0,0,1)} \circ
        \beta^{-1}_{(0,1,0)} \circ \beta_{(1,0,0)}$ sends
$1$ to $7$,  Proposition~\ref{p:ext} is violated.
Thus, this partial action does not extend to a partial
action of $\bbZ^3$.
Of course, if we replaced $\bbZ^3$ with the free group on
three generators, then the partial action would extend.
\end{example}

\begin{example} \label{ext}
We  give an example of an extension of a partial action from
 a positive cone which does not totally order the group.
This example will play a fundamental role in
Example~\ref{toroidal}, where we display Power's toroidal limit
 algebras as an analytic partial crossed product.

Let $X = \prod_1^i \{0,1\}$.  We could construct the same example on
any finite set with an even number of points.  For the purposes of
this example, we could as well denote the elements of the set with the
integers from 1 to $2^i$, but binary notation will facilitate the
discussion of the inverse system in example~\ref{toroidal}.  The group
in this example is $\bbZ^2$ and the positive cone is
$\Sigma =\{(a,b) \mid a \geq 0, b \geq 0\}$.  Write
$e_1 = (1,0)$ and $e_2 = (0,1)$.  Let $\omega$ be the finite odometer
map on $X$: $\omega(1,\dots,1)=(0,\dots 0)$ and otherwise $\omega$
adds 1 to the first entry with carries to the right.
Of course, $\omega$ is just a cyclic permutation of $X$.
Let 0 denote the identity element of the group, viz. $(0,0)$.
Also, let $D^0 = \{x \in X \mid x_1 = 0 \}$ and
$D^1 = \{x \in X \mid x_1 = 1 \}$.

Let $\beta(0)$ be the identity map on $X$; $\beta(e_1)$, the
restriction of $\omega$ to $D^0$ and $\beta(e_2)$, the restriction of
$\omega^{-1}$ to $D^0$.  For all other $s \in \bbZ^2$, 
$\beta(s)$ is empty.  In particular, we have
\begin{align*}
&D^0 \overset{\beta(e_1)}{\longrightarrow} D^1, \\
&D^0 \overset{\beta(e_2)}{\longrightarrow} D^1. 
\end{align*}
It is easy to check that $\beta$ is a partial action of 
$\Sigma$ on $X$.

We now wish to extend $\beta$ to an action of $\bbZ^2$ on $X$.  
Proposition~\ref{p:ext} gives guidance on how to proceed.
For example, both $\beta(e_1) \circ \beta(-e_2)$ and
$\beta(-e_2) \circ \beta(e_1)$ must be restrictions of
$\beta(e_1 - e_2)$.  But the first composition is
$\omega^2$ restricted to $D^0$ and the second is
$\omega^2$ restricted to $D^1$.  Thus, we have no choice but to take
$\beta(e_1 - e_2) = \omega^2$.  In a similar vein,
$\beta(e_1) \circ \beta(-e_2) \circ \beta(e_1)$ must be a restriction
of $\beta(2e_1 - e_2)$.  This is actually the only constraint, so we
will be able to take $\beta(2e_1 - e_2)$ equal to
$\omega^3$ restricted to $D^0$.  Considerations of this sort suggest
that the following will be an extension of $\beta$ from $\Sigma$ to
$\bbZ$:
\begin{align*}
&\beta(ae_1 - ae_2) = \omega^{2a}, \text{ for all $a \in \bbZ$,}\\
&\beta((a+1)e_1 - ae_2) \text{ is the restriction of 
$\omega^{2a+1}$ to $D^0$, for all $a \in \bbZ$,}\\
&\beta(ae_1 -(a+1)e_2) \text{ is the restriction of
$\omega^{2a+1}$ to $D^1$, for all $a \in \bbZ$,}\\
&\beta(s) \text{ is empty for all other $s \in \bbZ^2$.}
\end{align*}

It is a routine matter to check that the extended $\beta$ is a partial
action of $\bbZ^2$ on $X$.  Since we are dealing with a group,
it is enough to check conditions $1'$, 2 and 3; the only sets which
can be domains are $X$, $D^0$, $D^1$ and $\emptyset$ so the
calculations are easy.  Of course, the third arrow condition is also
easy to check.

If we visualize $\bbZ^2$ as the lattice points of a plane, and if we
put $\beta(ae_1+be_2)$ at node $(a,b)$,
 we can ``picture'' the partial action as follows:
\[
\begin{matrix}
\omega^{-4} & \omega^{-3}|_{D^0} & \emptyset & 
\emptyset & \emptyset \\[.25ex]
\omega^{-3}|_{D^1} & \omega^{-2} & \omega^{-1}|_{D^0}
     & \emptyset & \emptyset \\[.25ex]
\emptyset & \omega^{-1}|_{D^1} & id_X & \omega^1|_{D^0}
     & \emptyset \\[.25ex]
\emptyset & \emptyset & \omega^1|_{D^1} & \omega^2
      & \omega^3|_{D^0} \\[.25ex]
\emptyset & \emptyset &  \emptyset & \omega^3|_{D^1}
     & \omega^4
\end{matrix}
\]

\end{example}

\section{Analytic Partial Crossed Products} \label{s:apcp}

We shall define an analytic partial crossed
 product as a subalgebra of
a \cstar partial crossed product.
 This applies to any partial action of a
cone $\Sigma$ which has a unique extension to a partial action of the
whole group.  In particular, the definition is available whenever
$\Sigma$ totally orders $G$ (Proposition \ref{prop:toextension}).

A brief review of the definition of a partial crossed product
\cstar algebra will provide notation and terminology.
 If $\beta$ is a
partial action of an abelian group $G$ on a \cstar algebra $A$, we
first consider the set $\sP$ of all formal polynomials of the form
$\sum f_n U^n$, where $f_n \in D_n$, for each $n$ in some finite
subset of $G$.  Define on this set a multiplication and an
involution.  The multiplication is a ``twisted'' convolution product;
it is determined by specifying the product of two monomials:
\[
  f U^n g U^m = \alpha_n(\alpha_{-n}(f)g)U^{n+m}.
\]
Here, $\alpha_{-n}(f) \in D_{-n}$, so the product
$\alpha_{-n}(f)g \in D_{-n} \cap D_m$.  From the definition of
partial action, $\alpha_n(\alpha_{-n}(f)g) \in D_{n+m}$.
The involution is defined (on monomials) by
\[
(fU^n)^* = \alpha_{-n}(\overline{f}) U^{-n}.
\]
Thus, $\sP$ is a \star algebra.

Typically, we think of $U^m$ as the partial isometry implementing 
$\alpha_m$.
Notice, however, that if $\alpha$ does not have the composition 
property, then the product of the partial isometries implementing 
$\alpha_m$ and $\alpha_n$ may be a proper restriction of the 
partial isometry implementing $\alpha_{m+n}$.

Define a norm on $\sP$ by
\[
\left\| \sum f_n U^n \right\|_L = \sum \|f_n\|.
\]
A \cstar norm on $\sP$ is defined by
\[
\|x\| = \sup_{\pi} \|\pi(x)\|,
\]
where the supremum is taken over all representations which are
continuous with respect to the $L$-norm.
The partial crossed product $A \times_{\alpha} G$ is
the completion of
$\sP$ with respect to the \cstar norm. (If, instead of polynomials, we
had used all formal power series
$\sum_{n \in G} f_n U^n$ for which $\sum_n \|f_n\|$ converges, then we
would have obtained a Banach \star algebra whose enveloping
\cstar algebra is $A \times_{\alpha} G$.  This is what appears in
\cite{MR95g:46122,MR96i:46083}.)

\begin{defn}
The \emph{analytic crossed product} $A \times_{\alpha} \Sigma$ is
defined to be the closure in $A \times_{\alpha} G$ of
those polynomials
$\sum f_n U^n$ for which all $n$ are in $\Sigma$.
\end{defn}

Since we are restricting attention to the case where
$A=C_0(X)$ is an abelian \cstar algebra  we rephrase the
definition of the multiplication in terms of the partial action
$\beta$ of $\Sigma$ on $X$, where $\beta$
is associated with $\alpha$:
\[
f U^n g U^m = [(f \circ \beta_n)g]\circ \beta_n^{-1} U^{n+m}.
\]
For future reference, note that if $h$ denotes the coefficient function
$[(f \circ \beta_n)g]\circ \beta_n^{-1}$, then
\[
h(x) =
\begin{cases}
f(x) g(\beta_n^{-1}(x)), &\text{if $x \in X_n \cap X_{n+m}$},\\
0, &\text{otherwise.}
\end{cases}
\]
When $\beta$ is a partial action on $X$ and $\alpha$ is the dual
action on $C_0(X)$, we may write either
$C_0(X) \times_{\alpha} G$ or $C_0(X) \times_{\beta} G$ for the
partial crossed product.

\begin{remark}
In place of the \star algebra $\sP$ used in defining a partial crossed
product we can use a somewhat smaller algebra.  Let $\sP_c$ be the
subalgebra of $\sP$ consisting of all elements of $\sP$ whose
coefficients have compact support.  So, a polynomial
$\sum f_n U^n$ is in $\sP_c$ if, and only if, each
$f_n \in D_n$ and is compactly supported.

When $\sP$ is provided with the \cstar norm, $\sP_c$ is a dense
subalgebra.  This remark will be useful in Theorem~\ref{T:bicor}.
\end{remark}

We devote the remainder of this section to examples of analytic
partial crossed product algebras.
These include analytic limit algebras, instances of a general result
we prove in the next section, and also a family of non-triangular algebras.

The partial actions of the groups in 
Examples~\ref{eg:std} and
\ref{eg:refine} both yield \mbox{$2^{\infty}$ UHF} \cstar algebras.
The two analytic subalgebras are, respectively, the standard embedding
TAF algebra and the refinement embedding TAF algebra.
More generally, the standard $\bbZ$-analytic algebras
 of~\cite{MR94m:46113b,MR95f:46113}, the locally 
constant cocycle nest algebras
of~\cite{MR98d:47095}, and the order preserving 
algebras of~\cite{MR96k:46099}
are all analytic partial crossed product algebras.
The validity of these assertions follows directly from
Theorem~\ref{T:bicor}.

Before constructing our non-triangular example, 
we need a preliminary result.

\begin{prop} \label{P:indlim}
Fix an abelian group $G$.
\begin{enumerate}
\item Let $\beta$ and $\hat{\beta}$ be two
partial actions of $G$ on locally compact Hausdorff spaces $X$ and
$\hat{X}$.
As usual, suppose that $\beta_g$ has domain $D_{-g}$ and
$\hat{\beta}_g$ has domain $\hat{D}_{-g}$, for each $g$.
If $\phi \colon \hat{X} \to X$ is a continuous surjection so that,
for all $g \in G$ and all $x \in \hat{X}$,
$\phi(x) \in D_g$ if, and only if, $x \in \hat{D}_g$ and
\begin{equation} \label{E:coher}
 \beta_g(\phi(x))=\phi(\hat{\beta}_g(x)),
\end{equation}
then there is a continuous injection
$\Phi \colon C_0(X) \times_{\beta} G 
\to C_0(\hat{X}) \times_{\hat{\beta}} G$.
\newline When equation~\eqref{E:coher} holds, 
we will say that $\phi$ \emph{intertwines}
the two actions.
\item
Suppose that $(X_i,\phi_i)$ is an inverse system of locally compact
Hausdorff spaces in which each
 $\phi_i \colon X_{i+1} \to X_i$ is a continuous
surjection and is proper (in the sense that the inverse images of
compact sets are compact).
Assume further that there is a family of
partial actions $\beta_i$ on the $X_i$ so
that each $\phi_i$ intertwines $\beta_{i+1}$ and $\beta_i$.
Let $\Phi_i \colon C_0(X_i) \times_{\beta_i} G \to
C_0(X_{i+1}) \times_{\beta_{i+1}} G$ be the map induced by Part~(1).
If $X = \varprojlim(X_i,\phi_i)$, then there is a partial
 action $\beta$
of $G$ on $X$ so that $C_0(X) \times_\beta G  $ is isomorphic to
$\varinjlim \left(C_0(X_i) \times_{\beta_i} G, \Phi_i \right)$.
\item
Suppose further that $\Sigma$ is a positive cone in $G$.  Then
$\Phi_i$ maps $C_0(X_i) \times_{\beta_i} \Sigma$ into
$C_0(X_{i+1}) \times_{\beta_{i+1}} \Sigma$ and
$C_0(X) \times_{\beta} \Sigma = \varinjlim C_0(X_i) 
\times_{\beta_i} \Sigma$.
\end{enumerate}
\end{prop}

\begin{proof}
We prove (1) first.
As is well-known, the map 
$\tilde{\phi} \colon C_0(X) \to C_0(\hat{X})$ given
by $f \mapsto f \circ \phi$ is a continuous injection.
Denote the monomials of $C_0(X) \times_{\beta} G$ by $f U^m$
and those of $C_0(\hat{X}) \times_{\hat{\beta}} G$ by $f V^m$.
We define $\Phi$ on the monomials of $C_0(X) \times_{\beta} G$ by
$\Phi(f U^m) = \tilde{\phi}(f) V^m$ and then extend to polynomials.
Since $\tilde{\phi}$ is an injection, so is $\Phi$.

To see that $\Phi$ is a homomorphism, observe first that if
$\alpha$ and $\hat{\alpha}$ are the dual partial actions on
$C_0(X)$ and $C_0(\hat{X})$, 
then equation~\eqref{E:coher} implies that
$\tilde{\phi} \circ \alpha_m = \hat{\alpha}_m \circ \tilde{\phi}$.
Thus,
\begin{align*}
\Phi( f U^m ) \Phi( g U^n )
&= \tilde{\phi}(f) V^m \tilde{\phi}(g) V^n \\
&= \hat{\alpha}_m(\hat{\alpha}_{-m}(\tilde{\phi}(f)) \tilde{\phi}(g) )
    V^{m+n} \\
&= \hat{\alpha}_m(\tilde{\phi}(\alpha_{-m}(f)) 
   \tilde{\phi}(g) ) V^{m+n}\\
&= \hat{\alpha}_m(\tilde{\phi}(\alpha_{-m}(f) g) ) V^{m+n}\\
&= \tilde{\phi}(\alpha_m(\alpha_{-m}(f) g) ) V^{m+n}\\
&= \Phi(\alpha_m(\alpha_{-m}(f) g)  U^{m+n}) \\
&= \Phi(f U^m g U^{n}).
\end{align*}

It remains only to show that $\Phi$ is continuous on the polynomials
with the \cstar norms,
so that we may extend by continuity to the crossed product.
Clearly, $\Phi$ in contractive when the polynomials are equipped with
the $L$-norm, since $\tilde{\phi}$ is contractive on $C_0(X)$.
Composing $\Phi$ with an $L$-norm continuous representation of the
polynomials of $C_0(\hat{X}) \times_{\hat{\beta}} G$ gives an $L$-norm
continuous representation of the 
polynomials of $C_0(X) \times_{\beta} G$,
and so $\Phi$ is contractive with respect to the \cstar norms.

\medskip

Turning to part (2), we can identify $X$
 with the set of sequences
$\mathbf{x} = (x_1,x_2,\ldots)$, where each $x_i \in X_i$ and
$\phi_i(x_{i+1})=x_i$ for all $i$.
Of course, the topology is the relative product topology.
For each $g \in G$, let
$\dom \beta(g) = \{\mathbf{x} \mid x_i \in \dom \beta_i(g)
\text{ for all }i\}$
and define $\beta(g)\mathbf{x}$ to be 
$(\beta_1(g)x_1,\beta_2(g)x_2,\ldots)$.
The intertwining condition implies that
$\phi_i(\beta_{i+1}(g)x_{i+1})=\beta_i(g)x_i$ for all $i$,
and so $\beta(g)\mathbf{x} \in X$.

For $w \in X_i$, let $J(w) = \{ (y_j) \in X \mid y_i = w \}$.
Note that $y_j$ for $j \le i$ are determined by the condition
$y_i=w$.  If $x \in \dom \beta_i(g)$ and $\mathbf y \in J(w)$,
then $y_j \in \dom \beta_j(g)$, for all $j$, and
$\mathbf y \in \dom \beta(g)$.  Thus
$J(w) \subseteq \dom \beta(g)$ and
$\beta(g)(J(w)) = J(\beta_i(g)w)$.

Let $C_i$ be the subalgebra of $C_0(X) \times_\beta G$ generated
by $f U^m$ where $m \in G$ and $f \in C_0(X)$ is constant on the
sets $J(w)$, $w \in X_i$.
Then $C_i$ is isomorphic to $C_0(X_i) \times_{\beta_i} G$ and
$C_i \subseteq C_{i+1}$, for all $i$.
Under these isomorphisms and containments, the diagram
\[
\begin{diagram}
\node{C_i} \arrow{e} \arrow{s}
\node{C_{i+1}} \arrow{s}\\
\node{C_0(X_i) \times_{\beta_i} G} 
\arrow{e,t}{\Phi_i} 
\node{C_0(X_{i+1}) \times_{\beta_{i+1}} G}
\end{diagram}
\]
commutes.
To prove  that $C_0(X) \times_\beta G  $ is isomorphic to
$\varinjlim \left(C_0(X_i) \times_{\beta_i} G, \Phi_i \right)$
it is sufficient to show that
 the union of the $C_i$ is dense in $C_0(X) \times_\beta G$.

 To this end, let
\[
F_i = \{f \in C_0(X) \mid f \text{ is constant on $J(w)$, for
each $w \in X_i$}\}.  
\]
Then $F_1 \subseteq F_2 \subseteq F_3 \subseteq \dots$ and
each $F_i$ is a \star subalgebra of $C_0(X)$.  If we show
that the subalgebra $\bigcup F_i$ is dense in
$C_0(X)$, then the density of $\bigcup C_i$ in
$C_0(X) \times_{\beta} G$ follows.
To do this, all we need to show is that $\bigcup F_i$ separates points
of $X$.

Let $x$ and $y$ be two distinct points in $X$.  Then there is
an index $i$ such that $x_i \neq y_i$.  Therefore,
$J(x_i) \cap J(y_i) = \emptyset$.  Let $g$ be a continuous function
on $X_i$ such that $g(x_i) \neq g(y_i)$.  The composition of
$g$ with the canonical projection of $X$ onto $X_i$ is then an element
of $F_i$ which separates $x$ from $y$.  This completes the proof of
the second part of the Proposition.

The proof of the third part is trivial.
\end{proof}

Proposition~\ref{P:indlim} provides a convenient framework for
realizing limit algebras as analytic partial crossed products,
as the following examples illustrate.

\begin{example} \label{standard}
To obtain the $2^{\infty}$-UHF algebra and the 
$2^{\infty}$-standard TUHF subalgebra as direct limits in a partial
crossed product framework, 
let  $X_i = \prod_1^i \{0,1\}$ and let
$\beta_i$ be the finite odometer restricted to 
$X_i\backslash\{(1,\ldots,1)\}$.
(The odometer map sends a tuple $(x_1, \ldots, x_i)$ to
$(y_1, \ldots, y_i)$, where $1 + \sum x_k 2^{k-1} = \sum y_k 2^{k-1}
\pmod{2^i}$.)  The powers of $\beta_i$ give a partial action of
$\bbZ$ on $X_i$.  Let $\phi_i \colon X_{i+1} \to X_i$ be the map which
deletes the last entry of elements of $X_{i+1}$.  Then
$C(X_i) \times_{\beta_i} \bbZ \cong M_{2^n}$, 
$C(X_i) \times_{\beta_i} \bbZ^+ \cong T_{2^n}$, and the
embeddings $\Phi_i$ are the standard embedding maps.  Thus, we recover
the usual presentation of the standard TUHF algebra.
\end{example}

\begin{example} \label{BunceDeddens}
As a second easy application, this time 
with (full) crossed products, let
$\alpha_i$ be the odometer on $X_i = \prod_1^i \{0,1\}$
and let $\phi_i$ be the same as the Example~\ref{standard}.
Then each $C(X_i) \times_{\alpha_i} \bbZ$ is isomorphic to
$M_{2^i}(C(\bbT))$---$\bbT$ is unit circle---and, since the induced
map $\Phi_i$ is the twice-around embedding, the direct limit is
the $2^{\infty}$-Bunce-Deddens \cstar algebra.
The inverse limit action of the 
$\alpha_i$ on $X = \prod_1^{\infty} \{0,1\}$
is the odometer map $\alpha$ on $X$.
This gives the known result that $C(X) \times_{\alpha} \bbZ$ is
isomorphic to $\varinjlim (M_{2^i}(C(\bbT)),\Phi_i)$.
\end{example}

\begin{example} \label{toroidal}
We realize the toroidal limit algebras of~\cite{MR93d:46092} as
analytic partial crossed products.
Let $X_i = \prod_1^i \{0,1\}$ and once again let $\phi_i$ send
$(x_1,\ldots,x_i,x_{i+1})$ to $(x_1,\ldots,x_i)$.
For clarity, we confine ourselves to the $2^\infty$ case;
replacing each $\{0,1\}$ with $\{0,\ldots,n_i\}$, where
$2 \mid n_1$, gives the general construction.

For each $i$, let $\beta_i$ be the partial action of $\bbZ^2$
and the positive cone $\Sigma$ described in Example~\ref{ext}.
It is easy to check that the $\phi_i$ and the $\beta_i$ satisfy the
hypotheses of Proposition~\ref{P:indlim}.

For each $i$, $C(X_i) \times_{\beta_i} \Sigma$ is a $2^i$-cycle
algebra, that is, the subalgebra of $M_{2^i}$ spanned by the diagonal
matrix units, the matrix unit $e_{1,2^i}$, and the matrix units
$e_{j,k}$ where $j$ is odd, $k$ is even, and $|j-k|=1$.
As $\phi_i$ is the usual double cover embeddings of $X_{i+1}$ into
$X_i$, the induced map $\Phi_i$ wraps the $2^i$-cycle
algebra twice around the $2^{i+1}$-cycle algebra.
Thus, the direct system $(C(X_i) \times_{\beta_i} \Sigma, \Phi_i)$
is precisely the direct system given by Power~\cite[p.~51]{MR93d:46092}
and appealing to his Theorem~4.1, we have a generating subalgebra
of the $2^{\infty}$-Bunce-Deddens \cstar algebra.

We can show this directly.
Consider the element 
\[
y_i =\chi_{D^1} U^{e_1}+\chi_{D^0}U^{-e_2}
\]
in $C(X_i) \times_{\beta_i} \bbZ^2$.  ($D^0$ and $D^1$ are the domain
and range sets for $\beta_i(e_1)$ and $\beta_i(e_2)$, as in
Example~\ref{ext}.)  We first observe that $y_i$ is unitary.
The calculation uses the following
facts (in which $e$ is either $e_1$ or
$e_2$ and $\beta_i(-e) = \beta_i(e)^{-1}$):
\begin{align*}
\chi_{D^1} \circ \beta_i(e) &= \chi_{D^0}, &
\chi_{D^0} \circ \beta_i(-e) &= \chi_{D^1}, \\
\chi_{D^0} \circ \beta_i(e) &= 0, &
\chi_{D^1} \circ \beta_i(-e) &= 0.
\end{align*}
We then have
\[
y_i^* = \chi_{D^1} \circ \beta_i(e_1) U^{-e_1}
 + \chi_{D^0} \circ \beta_i(-e_2) U^{e_2}
= \chi_{D^0} U^{-e_1} + \chi_{D^1} U^{e_2}
\]
and
\begin{align*}
y_i^{\vphantom{*}} y_i^*
 &= (\chi_{D^1} U^{e_1}+\chi_{D^0}U^{-e_2})
  (\chi_{D^0} U^{-e_1} + \chi_{D^1} U^{e_2})\\
&= \chi_{D^1} \chi_{D^0} \circ \beta_i(-e_1) U^{e_1 -e_1}
+\chi_{D^1} \chi_{D^1} \circ \beta_i(-e_1) U^{e_1+e_2} \\
&\quad +\chi_{D^0} \chi_{D^0} \circ \beta_i(e_2) U^{-e_2-e_1}
+\chi_{D^0} \chi_{D^1} \circ \beta_i(e_2) U^{-e_2+e_2} \\
&= \chi_{D^1} U^0 + \chi_{D_0} U^0 \\
&= \chi_{X_i} U^0 = I.
\end{align*}
Similarly, $y_i^* y_i^{\vphantom{*}}=I$.

Observe that $y_i$ implements the action of the odometer
on $X_i$ (transferred to $C(X_i)U^0$): for $f \in C(X_i)$,
\begin{align*}
y_i^* fU^0 y_i^{\vphantom{*}} &=
(\chi_{D^0} U^{-e_1} + \chi_{D^1} U^{e_2}) fU^0
(\chi_{D^1} U^{e_1}+\chi_{D^0}U^{-e_2}) \\
&=\chi_{D^0} f \circ \beta_i(e_1) \chi_{D^0} U^0
  + \chi_{D^1} f \circ \beta_i (-e_2) \chi_{D^1} U^0 \\
&= \left[\chi_{D^0} f \circ \omega|_{D^0} +
\chi_{D^1} f \circ w|_{D^1}   \right]U^0 \\
&= f \circ \omega U^0.
\end{align*}

We next claim that $y_i$ and $C(X_i)U^0$ generate
$C(X_i) \times_{\beta_i} \bbZ^2$.  Let $A$ be the subalgebra generated
by $y_i$ and $C(X_i)U^0$.  Multiply $y_i$ on the left
by $\chi_{D^1}U^0$ to see that $\chi_{D^1}U^{e_1} \in A$ and on the
left by $\chi_{D^0}U^0$ to see that 
$\chi_{D^0} U^{-e_2} \in A$.  The
adjoint of the latter, $\chi_{D_1} U^{e_2}$ is therefore also in $A$. 
The square of $y_i$ is $\chi_{X_i} U^{e_1-e_2}$, so this monomial is
in $A$.  (We omit this and subsequent calculations, since they are
similar to the ones done above.)  For any positive integer $a$, the 
$a^{\text{th}}$ power of $\chi_{X_i} U^{e_1-e_2}$ is
$\chi_{X_i} U^{ae_1-ae_2}$; it follows that
 $\chi_{X_i} U^{ae_1-ae_2} \in A$ for all integers $a$.
Multiplication of $\chi_{X_i} U^{ae_1-ae_2}$ on the left
by $\chi_{D^1}U^{e_1}$  yields
$\chi_{D^1} U^{(a+1)e_1 - ae_2} \in A$; multiplication on the left 
by $\chi_{D^0} U^{-e_2}$ yields
$\chi_{D^0} U^{ae_1-(a+1)e_2} \in A$.  Now if $s$ is an element 
of $\bbZ^2$ which is not of the form $ae_1 - ae_2$ or
$(a+1)e_1 -ae_2$ or $ae_1 -(a+1)e_2$, then $\beta_i(s)$ is empty;
therefore, the only possible coefficient for $U^s$ in the partial
crossed product is 0.  Thus, we have shown that every monomial
$\chi_{\ran \beta_i (s)} U^s$ is in $A$.  But these monomials and
$C(X_i)U^0$ generate $C(X_i) \times_{\beta_i} \bbZ^2$, so
$A = C(X_i) \times_{\beta_i} \bbZ^2$.

Let $\Psi_i$ denote the inclusion of $C(X_i) \times_{\beta_i} \bbZ^2$
in $C(X_{i+1}) \times_{\beta_{i+1}} \bbZ^2$.
As $\Psi_i(y_i)=y_{i+1}$, we may let $y$ denote the image
in $C(X) \times_\beta \bbZ^2$ of (any) $y_i$.
Since $C(X_i)U^0$ and $y_i$ generate $C(X_i) \times_{\beta_i} \bbZ^2$
for each $i$, it follows that $C(X)U^0$ and $y$ generate
$C(X) \times_{\beta} \bbZ^2$.
Letting $\alpha$ be the action of the odometer map on $X$,
it follows that $C(X) \times_{\beta} \bbZ^2$ is a quotient of
$C(X) \times_{\alpha} \bbZ$, the $2^\infty$ Bunce-Deddens
\cstar algebra.
But the latter algebra is simple, so $C(X) \times_{\beta} \bbZ^2$
is the $2^\infty$ Bunce-Deddens algebra.  
By Proposition~\ref{P:indlim}, $C(X) \times_{\beta} \Sigma$ is Power's
toroidal limit algebra.
\end{example}

\section{Analyticity and Coordinates} \label{s:anco}

Theorem~\ref{T:bicor} relates partial crossed products to groupoid
\cstar algebras, with the analytic subalgebras in correspondence.
Relevant aspects of the groupoid \cstar algebra construction will
be reviewed briefly in the course of the proof.
For complete, systematic accounts, see Renault~\cite{MR82h:46075} or
Paterson~\cite{MR2001a:22003}.
For a convenient summary, see Muhly and Solel~\cite{MR90m:46098}.

Recall that a partial action is said to be \emph{free} if
$\alpha_t(x) = x$ implies that $t=0$.

\begin{thm} \label{T:bicor}
Let $G$ be a countable discrete subgroup of $\bbR$ with positive cone
$\Sigma = G \cap [0,\infty)$.  Let $\alpha$ be a free
partial action of $G$ on
a separable abelian \cstar algebra $A$.  Then there is a locally
compact $r$-discrete principal groupoid $\sG$ with second countable
unit space and a locally constant real valued cocycle on $\sG$ such
that the partial crossed product
$A \times_{\alpha} G$ is \star isomorphic to the groupoid
\cstar algebra $C^*(\sG)$; the analytic subalgebra of
$A \times_{\alpha} G$ associated with $\Sigma$ is carried by this
isomorphism onto the analytic subalgebra of $C^*(\sG)$ determined by
the cocycle.

Conversely, given a locally compact $r$-discrete principal groupoid
$\sG$ with second countable unit space and a locally constant
real valued cocycle, there is a countable discrete subgroup $G$ of
$\bbR$ and a partial action $\alpha$ of $G$ on a
 locally compact second countable Hausdorff space $X$
such that $C^*(\sG)$ is isomorphic to
$C_0(X) \times_{\alpha} G$; again, the analytic subalgebra of
$C^*(\sG)$ determined by the cocycle is carried by this isomorphism
onto the analytic subalgebra of
$C_0(X) \times_{\alpha} G$ determined by the positive cone
$G \cap [0,\infty)$.
\end{thm}

\begin{proof}
Let $A = C_0(X)$ be a separable
abelian \cstar algebra and let $G$ be a
discrete subgroup of $\bbR$ with a free partial action $\alpha$ on
$X$.
As a set, the groupoid $\sG$ is:
\[
\sG = \{ (x,\alpha_t(x)) \mid t\in G, x \in X_{-t} \}.
\]
This is an equivalence relation on $X$
(transitivity follows from the definition of partial action, most
trivially from the third arrow property),
 so the groupoid is principal.

For each $t \in G$ and each open subset
$U \subseteq X_{-t}$, let
$\sO_{t,U} = \{(x,\alpha_t(x)) \mid x \in U \}$.  Since the action is
free, $\alpha_t(x) = \alpha_s(x)$ implies $t=s$; consequently, the
family $\{\sO_{t,U}\}$ is closed under finite intersections.
(Either $\sO_{t,U} \cap \sO_{s,V} = \emptyset$ or $s=t$ and
$\sO_{t,U} \cap \sO_{t,V} = \sO_{t,U \cap V}$.) Give $\sG$ the smallest
topology in which all the $\sO_{t,U}$ are open sets.  The family of all
$\sO_{t,U}$ is a basis for this topology.

With this topology, $\sG$ is a locally compact topological space.
For each $t \in G$, the graph of $\alpha_t$, namely
$\{(x,\alpha_t(x)) \mid x \in X_{-t} \}$ is an open subset of $\sG$.
Furthermore, this set (with the relative topology)
is  isomorphic to $X_{-t}$ (and $X_t$).  Indeed, the range and domain
maps on the groupoid when restricted to the graph of $\alpha_t$ yield
isomorphisms with $X_{-t}$ and $X_t$, respectively.

If $g_n$ is a convergent sequence in $\sG$, then for some $t \in G$,
$g_n$ is eventually in the graph of $\alpha_t$.  So, for all large
$n$, $g_n = (x_n, \alpha_t(x_n))$ with $x_n \in X_{-t}$.  Since
$g_n$ is convergent, so is $x_n$.  Thus, there is $x \in X_{-t}$ such
that $x_n \to x$ and $\alpha_t(x_n) \to \alpha_t(x)$.

Since the set $\sO_{t,U}$ corresponds to the set
$\sO_{t^{-1},\alpha_t(U)}$ under the inverse map
$x \to x^{-1}$, the inverse map is continuous (indeed, a homeomorphism
of $\sG$ onto itself).

To show that $\sG$ is a topological groupoid it remains to verify that
the multiplication is continuous (when $\sG^2$ is given the relative
product topology from $\sG \times \sG$).  Suppose $g_n$ and $h_n$ are
two convergent sequences in $\sG$ and that, for each $n$, $g_n$ and
$h_n$ are composable.  It follows that (for all large $n$), there
exist $t,s \in \sG$ and $x_n \in X_{-t}$ such that
$g_n = (x_n, \alpha_t(x_n))$ and
$h_n = (\alpha_t(x_n), \alpha_s(\alpha_t(x_n)))$ and further, that
$x_n$ is convergent in $X_{-t}$, $\alpha_t(x_n)$ is convergent in
$X_t$, and $\alpha_s(\alpha_t(x_n))$ is convergent in $X_{s+t}$.
  Consequently,
$g_n h_n = (x_n, \alpha_{s+t}(x_n))$ is convergent in $\sG$.  Thus,
multiplication is continuous and $\sG$ is a topological groupoid.

The unit space $\sG^0$ for the groupoid is the graph of $\alpha_0$,
i.e., $\{(x,x) \mid x \in X\}$ and hence is open.  Thus $\sG$ is a
locally compact $r$-discrete principal groupoid.

Define a cocycle $c$ on $\sG$ by $c(x,\alpha_t(x)) = t$, for all
$(x,\alpha_t(x)) \in \sG$.  Since $\alpha$ is a partial action, $c$
satisfies the cocycle property.  Clearly, $c$ is constant on each open
set $\sO_{t,U}$, so $c$ is a locally constant  real valued cocycle.

We need to verify that the groupoid \cstar algebra constructed from
$\sG$ is isomorphic with the partial crossed product \cstar algebra
induced by the partial action $\alpha$.  The construction of
$C^*(\sG)$ begins with the family $C_c(\sG)$ of continuous functions
on $\sG$ with compact support.  Since the graphs of the $\alpha_t$
form a disjoint family of open sets, the support of a function in
$C_c(\sG)$ intersects only finitely many graphs.

The set $C_c(\sG)$ is provided with an involution, given by the
formula

\[
f^*(x,\alpha_t(x)) = \overline{f(\alpha_t(x),x)}
\]
and a multiplication defined by
\[
f\cdot g(x,\alpha_t(x)) = \sum f(x,\alpha_s(x))
 g(\alpha_s(x), \alpha_t(x)),
\]
where the sum is taken over all $s$ for which $x \in X_{-s}$ and
$\alpha_s(x) \in X_{s-t}$.  Since $f$ and $g$ have compact support,
only finitely many terms in the sum are non-zero.  Furthermore, when
$f$ and $g$ are supported on $\sG$-sets, at most one term in this sum
is non-zero.

$C_c(\sG)$ is provided with a norm
(usually called the $I$-norm), in which the norm of $f$
is given by
\[
\|f\|_I = \max \left\{\sup_{x \in X} \sum_t |f(x,\alpha_t(x))|,
\,\sup_{y \in X} \sum_t |f(\alpha_t(y),y)|\right\}
\]
and then a \cstar norm is obtained by defining
$\|f\| = \sup_{\pi} \|\pi(f)\|$, where $\pi$ varies over all \star
representations of $C_c(\sG)$ which are norm decreasing with respect
to the $I$-norm.

Recall that the analytic subalgebra associated with the cocycle
$c$ is the closure in the \cstar norm of all the functions in
$C_c(\sG)$ which are supported on
$\{(x,y) \mid c(x,y) \geq 0 \}$.  For a discussion of
analytic subalgebras, see~\cite{MR90m:46098}.

We next  define a
\star isomorphism $\Phi$ between $\sP_c$
and $C_c(\sG)$. As a map from $(\sP_c, \|\ \|_L)$ to
$(C_c(\sG), \|\ \|_I)$, $\Phi$ will be norm decreasing but not norm
preserving.  But with respect to the \cstar norms on $\sP_c$ and
$C_c(\sG)$, $\Phi$ will be an isometry; its
 extension to the completions
yields a \star isomorphism from $A \times_{\alpha} G$ onto $C^*(\sG)$.

The isomorphism (and its properties) is determined by its action
on monomials in $\sP_c$.
So, if $fU^n$ is a monomial in $\sP_c$, define $\Phi(fU^n)$ to be the
function on $\sG$ given by
\[
\Phi(fU^n)(x,\alpha_t(x)) =
\begin{cases}
f(x), &\text{if $t=-n$ and $x \in X_n$,} \\
0, &\text{otherwise.}
\end{cases}
\]
Since $f$ is continuous and has compact support in $X_n$,
$\Phi(fU^n) \in C_c(\sG)$ and has support on the graph of $\alpha_{-n}$.

Since $(fU^n)^* = \overline{f\circ \alpha_n} U^{-n}$, we have
\begin{align*}
\Phi((fU^n)^*)(x,\alpha_t(x))
&= \Phi(\overline{f\circ \alpha_n} U^{-n})(x,\alpha_t(x)) \\
&=
\begin{cases}
\overline{f(\alpha_n(x))}, &\text{if $t=n$ and
$x \in X_{-n}$,} \\
0, &\text{otherwise,}
\end{cases}
\end{align*}
while
\begin{align*}
\Phi(fU^n)^*(x,\alpha_t(x))
&= \overline{\Phi(fU^n)(\alpha_t(x),x)} \\
&=
\begin{cases}
\overline{f(\alpha_n(x))}, &\text{if $t=n$ and $x \in X_{-n}$,} \\
0, &\text{otherwise.}
\end{cases}
\end{align*}
Thus $\Phi$ is \star preserving on monomials; when $\Phi$ is extended
to $\sP_c$ by linearity, it remains \star preserving.

Now suppose that $fU^n$ and $gU^m$ are monomials in $\sP_c$.  Recall
that $fU^ngU^m = hU^{n+m}$, where
\[
h(x) =
\begin{cases}
f(x)g(\alpha_{-n}(x)), &\text{if $x \in X_n \cap X_{n+m}$},\\
0, &\text{otherwise.}
\end{cases}
\]
Therefore
\begin{align*}
\Phi(fU^ngU^m)(x,\alpha_t(x))
&=\Phi(hU^{n+m})(x,\alpha_t(x))\\
&=
\begin{cases}
h(x), &\text{if $t=-n-m$ and $x \in X_{n+m}$,}\\
0, &\text{otherwise,}
\end{cases}\\
&=
\begin{cases}
f(x)g(\alpha_{-n}(x)),
 &\text{if $t=-n-m$ and $x \in X_n \cap X_{n+m}$,}\\
0, &\text{otherwise,}
\end{cases}
\end{align*}
On the other hand, since
$\Phi(fU^n)$ is supported on the graph of $\alpha_{-n}$ and
$\Phi(gU^m)$ is supported on the graph of $\alpha_{-m}$,
$\Phi(gU^m)\cdot\Phi(fU^n)(x,\alpha_t(x)) = 0$ whenever $t \neq -n-m$
and
\begin{align*}
\Phi(fU^n)\cdot\Phi(gU^m)(x,\alpha_{-n-m}(x))
&= \Phi(fU^n)(x,\alpha_{-n}(x))
\Phi(gU^m)(\alpha_{-n}(x),\alpha_{-m}(\alpha_{-n}(x))) \\
&=
\begin{cases}
f(x)g(\alpha_{-n}(x)),
&\text{if  $x \in X_n \cap X_{n+m}$,}\\
0, &\text{otherwise.}
\end{cases}
\end{align*}
Thus $\Phi(fU^ngU^m)=\Phi(fU^n)\Phi(gU^m)$.  It  follows
immediately  that $\Phi(pq) = \Phi(p)\Phi(q)$ for all
$p,q \in \sP_c$ and so
$\Phi$ is a \star isomorphism of $\sP_c$ onto $C_c(\sG)$.
 (To see that $\Phi$
 is surjective, recall that an element $f$ in $C_c(\sG)$ is
supported on only finitely many graphs; restrict $f$ to each of
these graphs.
 On each graph, the map $(x,\alpha_{-n}(x)) \mapsto x$ is a
 homeomorphism, so $f$ can be transferred to a continuous function with
 compact support on $X_n$.  These are the coefficients of the
 polynomial in $\sP_c$ which is mapped by $\Phi$ to $f$.)

Let $p = \sum f_nU^n$ be a polynomial in $\sP_c$.  Each
$\Phi(f_nU^n)$ is supported on the graph of a single $\alpha_t$ (viz.,
$t=-n$), so

\[
\Phi(p)(x,\alpha_t(x)) =
\begin{cases}
f_{-t}(x), &\text{if $-t$ is an index in the sum for $p$},\\
0, &\text{otherwise.}
\end{cases}
\]
Therefore, for any $x$,
\[
\sum_t|\Phi(p)(x,\alpha_t(x))| =
\sum_n|f_n(x)| \leq \sum_n\|f_n\| = \|p\|_L.
\]
(The first sum is taken over those $t$ for which $x \in X_{-t}$.  The
second sum is similarly restricted, but the third sum is taken over
all the indices in the expression for $p$.)  Similarly,
$\sum_t|\Phi(p)(\alpha_t(y),y)| \leq \|p\|_L$, for each $y$; it
follows that $\|\Phi(p)\|_I \leq \|p\|_L$.

Although the map $\Phi$ is not norm preserving, as we show after
the proof, the fact that
it does satisfy $\|\Phi(p)\|_I \leq \|p\|_L$ for all
$p \in \sP_c$
 is enough to imply that $\Phi$ is norm decreasing with respect
to the \cstar norms on $\sP_C$ and $C_c(\sG)$.
Indeed, if $\rho$ is any $\|\ \|_I$-decreasing representation of
$C_c(\sG)$, then $\rho \circ \Phi$ is a $\|\ \|_L$-decreasing
representation of $\sP_c$.
It follows that $\rho \circ \Phi$ is norm decreasing with respect to
the \cstar norm on $\sP_c$: $\|\rho(\Phi(p))\| \leq \|p\|$.
Since this is true for all $\|\ \|_I$-decreasing representations,
$\|\Phi(p)\| \leq \|p\|$.

We can, in fact,  show that $\Phi$ is an isometry with respect to
the \cstar norms on $\sP_c$ and $C_c(\sG)$.  This task is made simpler
by the fact that $G$ is an amenable group.  For in this case, the
reduced partial crossed product (as defined by McClanahan in
\cite{MR96i:46083}) is isomorphic to $A \times_{\alpha} G$.
This means that
we do not need to use all $\|\ \|_L$-decreasing representations
of $\sP_c$ to determine
the \cstar norm on $A \times_{\alpha} G$ .  Instead, we can
restrict to a family of
representations constructed from representations of $A$ acting on a
Hilbert space $H$ and the left regular representation of the group $G$
acting on $\ell^2(G,\sH)$.  (In fact, it will suffice to consider only
certain of these representations.)

Given a representation $\pi \colon A \to \sB(\sH)$, McClanahan
constructs a ``regular representation'' $\tilde{\pi}$ of $A$ acting on
$\ell^2(G,\sH)$.  This is done as follows.  For each $g \in G$, a
representation $\pi_g \colon D_g \to \sB(\sH)$ is defined by
$\pi_g(z) = \pi(\alpha_{-g}(z))$.  This is extended to a
representation of all of $A$ (in a unique way) by taking
$\pi'_g(z) = \text{s-}\lim_{\nu} \pi_g(u_{\nu}z)$, where
$u_{\nu}$ is an approximate identity in $D_g$ and
$\text{s-}\lim$ refers to the limit in the strong operator
topology.  Finally,
$\tilde{\pi}$ is defined by
$\tilde{\pi}(z) \xi(g) = \pi'_g(z)\xi(g)$ for all
$\xi \in \ell^2(G,\sH)$.  Let $\lambda$ denote the left regular
representation of $G$ acting on $\ell^2(G,\sH)$.  Then the
representations which determine the reduced partial crossed
product norm
(and hence, in our case, the partial crossed
 product norm) are those of
the form $\tilde{\pi} \times \lambda$.  The action of\
$\tilde{\pi} \times \lambda$ is given on monomials by
$\tilde{\pi} \times \lambda (f_tU^t) = \tilde{\pi}(f_t)\lambda_t$.

The situation is further simplified by Proposition 3.4 in
\cite{MR96i:46083}, which tells us that if $\pi$ is a faithful
representation of $A$, then $\tilde{\pi} \times \lambda$ is a faithful
representation of $A \times_\alpha G$.  We can obtain a faithful
representation of $A$ from the following family of representations:
for each $x\in X$, let $\pi^x \colon A = C_0(X) \to \bbC$ be defined
by $\pi^x(f) = f(x)$; just let
$\pi = \sum^{\oplus}\pi^x$.  McClanahan's construction respects direct
sums, so
$\tilde{\pi} \times \lambda = \sum^{\oplus}\widetilde{\pi^x} \times
\lambda$.

In order to prove that $\Phi$ is norm decreasing (with respect to the
\cstar norms), it will suffice to show that for $p \in \sP_c$,
$\|(\tilde{\pi}\times\lambda)(p)\| \leq \|\Phi(p)\|_I$.
 From this, it follows
that $(\tilde{\pi}\times\lambda) \circ \Phi^{-1}$
is continuous with respect to
the $\|\ \|_I$-norm, and hence with respect to the \cstar norm.  But
$\tilde{\pi}\times \lambda$ is a faithful
representation, so $\Phi^{-1}$ is
\cstar norm decreasing; since the same is true for $\Phi$, $\Phi$ is
an isometry from $\sP_c$ onto $C_c(\sG)$.  Therefore, $\Phi$ has a
unique extension to a \star isomorphism from $A \times_{\alpha} G$ to
$C^*(\sG)$.  This extension clearly maps the analytic algebra
associated with the positive cone in $G$ onto the analytic algebra
associated with the cocycle on $\sG$.

Since
$\tilde{\pi} \times \lambda = \sum^{\oplus}\widetilde{\pi^x} \times
\lambda$, we can complete this direction of the proof by showing
that, for any $x \in X$ and $p \in \sP_c$,
$\|(\widetilde{\pi^x} \times \lambda) (p)\| \leq \|\Phi(p)\|_I$.
 Now, for each
$t \in G$, $\pi^x_t \colon C_0(X_t) \to \sB(\bbC) \cong \bbC$ is given
by
\[
\pi^x_t(f) =
\begin{cases}
\alpha_{-t}(f)(x) =  f(\alpha_t(x)), &\text{if } x \in X_{-t}, \\
0, &\text{otherwise}.
\end{cases}
\]

Hereafter,  $f(\alpha_t(x))$ will be understood to designate $0$
whenever $x \notin X_{-t}$.  The unique extension of $\pi^x_t$ from
$C_0(X_t)$ to all of $A=C_0(X)$ is given by exactly the same formula.
The regular representation $\widetilde{\pi^x}$ determined by $\pi^x$
maps $A$ into $\sB(\ell^2(G,\bbC))$ and is given by the formula
$(\widetilde{\pi^x}(f)\xi)(t) = \pi^x_t(f)\xi(t) =
f(\alpha_t(x))\xi(t)$.  We can more conveniently describe
$\widetilde{\pi^x}$ in terms of its action on the canonical basis
$\{\delta_t\}$ for $\ell^2(G,\bbC)$:
\[
\widetilde{\pi^x}(f)\delta_t = f(\alpha_t(x))\delta_t.
\]
The left regular representation $\lambda$, on the other hand, is given
by
\[
\lambda_t\delta_s = \delta_{t+s}.
\]
Thus, the matrix which represents  $\widetilde{\pi^x}$ with
respect to the canonical basis is a diagonal matrix and the
matrix which represents $\lambda_t$ is supported on a single sub or
super diagonal (main diagonal when $t=0$), where all the entries are
$1$'s.

If $f_tU^t$ is a monomial in $\sP_c$, then
\[
(\widetilde{\pi^x}\times\lambda)(f_tU^t)=
\widetilde{\pi^x}(f_t)\lambda_t.
\]
This is an operator in $\sB(\ell^2(G,\bbC))$ whose matrix has non-zero
entries only in the `diagonal' determined by $t$; in fact,
\begin{align*}
(\widetilde{\pi^x}\times\lambda) (f_tU^t)(\delta_s) &=
\widetilde{\pi^x}(f_t)\lambda_t \delta_s \\
&= \widetilde{\pi^x}(f_t)\delta_{t+s} \\
&= f_t(\alpha_{t+s}(x))\delta_{t+s}.
\end{align*}
If we let $k \in G$, the single entry in the $k$-column of
$\widetilde{\pi^x}(f_t)\lambda_t$ is
$f_t(\alpha_{t+k}(x))$ and the single entry in the $k$-row is
$f_t(\alpha_k(x))$.

Let $p=f_{t_1}U^{t_1}+\dots+f_{t_n}U^{t_n}$ be an arbitrary polynomial
in $\sP_c$.  Then
\[
(\widetilde{\pi^x} \times \lambda)(p) =
\widetilde{\pi^x}(f_{t_1}) \lambda_{t_1} + \dots +
\widetilde{\pi^x}(f_{t_n}) \lambda_{t_n}
\]
has finitely many entries in each row and in each column.
The $\ell^1$-norm of the $k$-column is
\[
c_k = |f_{t_1}(\alpha_{t_1+k}(x))| + \dots
+ |f_{t_n}(\alpha_{t_n+k}(x))|
\]
and the $\ell^1$-norm of the $k$-row is
\[
r_k = |f_{t_1}(\alpha_k(x))| + \dots + |f_{t_n}(\alpha_k(x))|.
\]
Therefore,
$\|(\widetilde{\pi^x}\times\lambda)(p)\| \leq \sup_k\{c_k,r_k\}$.  We
can complete the argument by showing that
$\sup_c\{c_k,r_k\} \leq \|\Phi(p)\|_I$.

Since $\Phi(p)$ is supported on the graphs of
$\alpha_{t_1}, \dots \alpha_{t_n}$ and since each term
$\Phi(f_{t_k}U^{t_k})$ is supported on the graph of
$\alpha_{t_k}$ alone,
\begin{align*}
\sum_t|\Phi(p)(z,\alpha_t(z))|
&=|\Phi(p)(z,\alpha_{-t_1}(z))| + \dots +
|\Phi(p)(z,\alpha_{-t_n}(z))|\\
&=|f_{t_1}(z)| + \dots + |f_{t_n}(z)|.
\end{align*}
It follows that
\[
\sup\{r_k\} \leq \sup_{z \in X} \sum_t
|\Phi(p)(z,\alpha_t(z))|.
\]
Also,
\begin{align*}
\sum_t|\Phi(p)(\alpha_t(z),z)|
&=\sum_t|\Phi(p)(\alpha_t(z),\alpha_{-t}(\alpha_t(z)))|\\
&=|f_{t_1}(\alpha_{t_1}(z))|+\dots+|f_{t_n}(\alpha_{t_n}(z))|.
\end{align*}
Hence,
\[
\sup\{c_k\} \leq \sup_{z\in X}\sum_t
|\Phi(p)(z,\alpha_t(z))|.
\]
This yields
$\|(\widetilde{\pi^x}\times\lambda)(p)\|\leq\|\Phi(p)\|_I$ and
the proof
of this direction of the theorem is complete.

The converse direction in the theorem
remains to be considered.  Start with an
$r$-discrete principal groupoid $\sG$ (based on a locally compact,
second countable, Hausdorff
space $X$) with a locally constant cocycle $c$.
Let $G$ be the range of
the cocycle.  Since $X$ is second countable and $c$ satisfies the
cocycle property, $G$ is a countable subgroup of $\bbR$.

For $t\in G$, the set
$U_t = \{(x,y) \mid c(x,y)=t\}$ is an open $\sG$-set in $\sG$.  With
$\pi_1$ and $\pi_2$ the natural projections of $X \times X$ onto $X$
(first and second coordinate projections), let $X_{-t} = \pi_1(U_t)$
and $X_t = \pi_2(U_t)$.  Then $U_t$ is the graph of a homeomorphism,
which we denote by $\alpha_t$, of $X_{-t}$ onto $X_t$.  The system
$\alpha$ satisfies all the requirements for a free partial action of
$G$ on $X$.  (Freeness follows from the fact that
$U_t \cap U_s = \emptyset$ whenever $t \neq s$.)

Since the groupoid $\sG$ is exactly the groupoid obtained from the
partial action $\alpha$ in the first part of the proof, the
groupoid \cstar algebra and the partial crossed product
\cstar algebra are isomorphic
(with appropriate correspondence of analytic subalgebras).
\end{proof}

\begin{remark}
\textbf{(1)}
The map $\Phi$ is not norm preserving with respect to the $I$-norm
and the $L$-norm, as the following example shows.
Let $X=\{1,\dots,n\}$.  Let $\beta_1$
be the map $k \mapsto k+1$ on the natural domain; $\beta_1$
generates a partial action of $\bbZ$ on $X$.  The partial crossed
product algebra (indeed, $\sP_c$) is isomorphic to $M_n$.  The
associated groupoid is $X \times X$, the full equivalence relation and
the groupoid \cstar algebra is again isomorphic to $M_n$.
If we identify both the partial crossed
 product algebra and the groupoid
\cstar algebra with $M_n$ then we can describe the $I$-norm and the
$L$-norm as follows:  for $a=(a_{ij})$
\begin{align*}
\|a\|_I &= \max(\max_i \sum_j |a_{ij}|, \max_j \sum_i |a_{ij}|), \\
\|a\|_L &= \sum_{t=-n+1}^{n-1}\max_{j-i=t}|a_{ij}|.
\end{align*}
In other words, to compute the $I$-norm, compute the $\ell^1$-norm of
each row and of each column in $a$ and take the largest value.  To
compute the $L$-norm, compute the $\ell^{\infty}$-norm of each diagonal
(determined by fixing values for $j-i$) and add all these numbers.

If we now take a matrix which is zero except for entries on the
counter diagonal, we obtain different values for the two norms.  To be
specific, suppose
\[
a_{ij} =
\begin{cases}
1, &\text{if $i+j=n+1$}, \\
0, &\text{otherwise}.
\end{cases}
\]
  Then
$\|a\|_I = 1$ and $\|a\|_L=n$.  This disparity
makes it clear that for
infinite dimensional algebras, the two norms
need not be equivalent.

\medskip
\noindent
\textbf{(2)}
The beginning of the proof of Theorem \ref{T:bicor} outlines the
passage from a partial action of a discrete abelian group on a locally
compact Hausdorff topological space to a groupoid \cstar algebra.
(The groupoid is the union of the graphs of the partial
homeomorphisms.) This connection has been studied for
$\bbZ$-partial actions by Peters and Poon in \cite{jpytp} and for
partial actions by countable discrete abelian groups by Peters and
Zerr in \cite{jprz}.  These papers use the language of partial
dynamical systems rather than partial actions, but there is little
substantive difference.  In particular, the question of when the
associated groupoid \cstar algebra is AF is addressed.  Peters and
Zerr give necessary and sufficient conditions on a partial dynamical
system for the associated groupoid \cstar algebra to be AF.  These
conditions are rather involved, so we won't restate them here, but
merely refer the reader to \cite{jprz}; in view of Theorem
\ref{T:bicor} these conditions characterize when a partial crossed
product by a discrete countable subgroup of $\bbR$ is an AF
\cstar algebra.

\medskip
\noindent
\textbf{(2)}
The order preserving normalizer of a TAF algebra plays an important
role in the study of the ideal structure of TAF algebras and in the
study of automatic continuity for algebraic isomorphisms
\cite{MR97e:47068,MR2000d:47103,MR2001k:47103}.
In particular, it is useful to know when an algebra is generated by
its order preserving normalizer.
The earliest result about when the order preserving normalizer
generates an algebra appears in \cite{MR92m:46089}, where it is shown
that if $\sG$ is an $r$-discrete principal groupoid with a continuous
cocycle onto a discrete ordered group, if $A$ is the analytic subalgebra
of $C^*(\sG)$ associated with the cocycle, and if $P \subset \sG$ is
the spectrum of $A$, then $P$ is the union of the monotone $\sG$-sets
which it contains.
Now, the (compact, open) $\sG$-sets correspond to the partial isometries
in $C^*(\sG)$ which normalize the diagonal $A \cap A^*$ and a $\sG$-set
is monotone if, and only if, the normalizing partial isometry is order
preserving.
Furthermore, the monotone $\sG$ sets cover $P$ if, and only if, $A$ is
generated by its order preserving normalizer.
Since locally constant cocycles are necessarily continuous,
the following corollary to Theorem~\ref{T:bicor} is immediate.

\begin{cor} \label{C:op}
If a TAF algebra is an analytic partial crossed product, then it is
generated by its order preserving normalizer.
\end{cor}
\end{remark}

\section{Conjugacy Results} \label{s:conjugacy}

A major theme in the study of analytic crossed products is that
$C(X) \times_\alpha \bbZ^+$ contains all of the
dynamical information about $\alpha \colon X \to X$, in the sense
that two homeomorphisms of $X$ are conjugate if, and only if,
the associated analytic crossed products are isomorphic.
See, for example,
\cite{MR35:1751,MR40:3322,MR86e:46063,MR89e:47069,%
MR90a:46175,MR93h:46095}.
In this section, we show that an analytic partial crossed product
$C(X) \times_\alpha \bbZ^+$ similarly contains all of the dynamical
information about the partial action $\alpha$.

There are several related results for limit algebras.
Power, in~\cite[Theorem~3.4]{MR93d:46092}, used inverse systems of
simplicial complexes to construct both operator algebras and
dynamical systems; he showed that conjugacy of the dynamical
systems is equivalent to isometric isomorphism of the operator
algebras and to isomorphism of the associated coordinates.
Poon and Wagner, in~\cite[Theorem~4.1]{MR95f:46113}, show that for
a family of subalgebras of crossed products with a distinguished point,
the algebras are isomorphic exactly when there is a conjugacy of
the dynamical systems that sends one distinguished point to the other.
Precisely, their systems are obtained from an essentially minimal
homeomorphism, i.e., one possessing unique minimal closed invariant set,
acting on a Cantor set; fixing a 
point $x$ in this unique minimal set, they
consider the subalgebra generated by the diagonal, $C(X)$, and the set
$\{ U f \mid f(x)=0, f \in C(X) \}$, where $U$ is canonical unitary
implementing the homeomorphism.

The composition of two partial homeomorphisms of a space $X$,
say $\alpha \colon A \to B$ and $\beta \colon C \to D$, is
the partial homeomorphism
$\alpha \circ \beta$ on the domain where this makes sense, i.e.,
on $C \cap \beta^{-1}(A)$.

\begin{defn}
We say two partial homeomorphisms
$\alpha \colon A \to B$ and $\beta \colon C \to D$ with
$A,B \subseteq X$ and $C,D \subseteq Y$
 are \emph{conjugate}
if there is a homeomorphism $\tau \colon X \to Y$ such that
$\tau$ maps the domain of $\alpha$ onto the domain of $\beta$, the
range of $\alpha$ onto the range of $\beta$ and
$\beta \circ \tau = \tau \circ \alpha $,
as partial homeomorphisms.

Two partial actions, $\alpha$ and $\beta$, of a
positive cone $\Sigma$ acting on
spaces $X$ and $Y$ are \emph{conjugate} if there
is a homeomorphism
$\tau \colon X \to Y$ such that $\tau$ induces a conjugacy between
$\alpha_s$ and $\beta_s$,
 for all $s \in \Sigma$.
\end{defn}

\begin{thm} \label{conjugacy}
Let $\alpha$ and $\beta$ be two partial actions of $\bbZ^+$  on
locally compact Hausdorff spaces $X$ and $Y$, respectively.  Assume
that each action satisfies the domain ordering property and that
each action is free.  Then
there is a continuous isomorphism of $C_0(X) \times_{\alpha} \bbZ^+$
onto $C_0(Y) \times_{\beta} \bbZ^+$ if, and only if, $\alpha$ and
$\beta$ are conjugate.
\end{thm}

\begin{proof}
The proof is a modification of the proof of Theorem~1 in
\cite{MR93h:46095}.  First consider the algebra
$A = C_0(X) \times_{\alpha} \bbZ^+$ associated with the action
$\alpha$.  Following the notation set up in Section~$\ref{s:anco}$,
$A$ is the closure of polynomials
$f_0 + f_1U^1+f_2U^2 + \dots + f_nU^n$, where each
$f_k \in C_0(X_k)$.  Let $\sC_1$ be the closed commutator ideal in
$A$ and,
for each $k \geq 2$, let $\sC_k$ be the closed ideal generated by
$k$-fold products of elements of $\sC_1$.
For $k \geq 1$, let $\sB_k$ be the closed ideal generated by all
polynomials of the form
$f_kU^k + \dots f_nU^n$, $n \geq k$.  Our first task is to show that
$\sC_k = \sB_k$, for all $k$.

Let $f_nU^n$ and $g_mU^m$ be two monomials in $\sP$.  Then
\[
f_nU^ng_mU^m = [(f_n \circ \alpha_n)g_m]\circ \alpha_{-n}U^{n+m}
\]
 and
\[
g_mU^mf_nU^n = [(g_m \circ \alpha_m)f_n]\circ \alpha_{-m}U^{n+m}.
\]
(The coefficients are supported in $X_n \cap X_{n+m}$
and $X_m \cap X_{n+m}$, but each of these sets is just
$X_{n+m}$ since $\alpha$ satisfies the domain ordering property.)

If $n=m=0$, then $f_0g_0 - g_0f_0 = 0$, since $C_0(X)$ is abelian.  If
either $m \neq 0$ or $n \neq 0$, then the commutator
$[f_nU^n, g_mU^m]$ lies in $\sB_1$.  By linearity, this is true for
all elements of $\sP$; hence $\sC_1 \subseteq \sB_1$.

If $f_0 \in C_0(X)$ and $g_m \in C_0(X_m)$ with $m \geq 1$, then
\[
[f_0, g_mU^m] =
[f_0g_m - [(g_m \circ \alpha_m)f_0]\circ \alpha_{-m}]U^m.
\]
Letting $h = f_0g_m - [(g_m \circ \alpha_m)f_0]\circ \alpha_{-m}$, we
have $h(x) = 0$ when $x \notin X_m$ and
\begin{align*}
h(x) &= f_0(x)g_m(x) - g_m(x) f_0(\alpha_{-m}(x)) \\
&= [f_0(x) - f_0(\alpha_{-m}(x))]g_m(x)
\end{align*}
when $x \in X_m$.

Since $\alpha_{-m}(x) \neq x$, for all $x$, we can, for each pair of
distinct points $x,y \in X_m$, choose $f_0 \in C_0(X)$ and
$g_m \in C_0(X_m)$ so that $h(y) \neq h(x) \neq 0$.
Let $\sE = \{g \in C_0(X_m) \mid gU^m \in \sC_1\}$.  Then $\sE$ is a
closed \star subalgebra of $C_0(X_m)$ and, from
what we have just observed,
$\sE$ separates points of $X_m$ and separates points from $\infty$.
Hence, by the Stone-Weierstrass theorem, $\sE = C_0(X_m)$.  As
$m \geq 1$ is arbitrary, $\sC_1 = \sB_1$.

If $p$ and $q$ are two polynomials in $\sC_1=\sB_1$,
then $pq$ is a
polynomial in $\sB_2$; thus $\sC_2 \subseteq \sB_2$.

Let $f_1U^1$ and $g_mU^m$ be two monomials (with $m \geq 1$).  Then
\[
f_1U^1g_mU^m = [(f_1 \circ \alpha_1)g_m]\circ \alpha_{-1}U^{m+1}.
\]
Letting $h = [(f_1 \circ \alpha_1)g_m]\circ \alpha_{-1}$, we have
$h(x) = 0$ when $x \notin X_{m+1}$ and
$h(x) = f_1(x)g_m(\alpha_{-1}(x))$ when
$x \in X_1 \cap X_{m+1} = X_{m+1}$.  (Since $\alpha$
satisfies the domain ordering property, $X_{m+1} \subseteq X_1$.)

The freedom to choose $f_1$ and $g_m$ arbitrarily in
$C_0(X_1)$ and $C_0(X_m)$ means that the set of
all coefficient functions
$h$ which arise in this fashion
separates points of $X_m$ and separates points from
$\infty$.  Again, the Stone-Weierstrass theorem shows that $\sC_2$
contains all polynomials of the form
$f_2U^2+\dots+f_nU^n$; hence $\sB_2 \subseteq \sC_2$.

Once we have shown that $\sC_{k-1} = \sB_{k-1}$, the same argument
yields $\sC_k = \sB_k$ (since $\sC_k$ is the closed ideal generated by
products of elements in $\sC_1$ and $\sC_{k-1}$); thus
$\sC_k = \sB_k$ for all $k \geq 1$.

 From the characterization of $\sC_1$ given above, it is evident
 that $A/\sC_1 \cong C_0(X)$.

Let $k \geq 1$.  We need to describe the ideals which lie between
$\sC_{k+1}$ and $\sC_k$ and are maximal in this class.  For each
$x \in X_k$, let
\[
I_x = \{f \in C_0(X) \mid f(x)=0 \text{ and } \supp f \subseteq X_k \}
\]
and for each $y \in X_{-k}$, let
\[
J_y = \{f \in C_0(X) \mid f(y)=0 \text{ and }
\supp f \subseteq X_{-k} \}
\]

Consider $I_x \sC_k + \sC_{k+1}$, the closure of all elements of the
form $fp+q$, where $f \in I_x$, $p \in \sC_k$, and $q \in \sC_{k+1}$.
Clearly,
$\sC_{k+1} \subseteq I_x\sC_k + \sC_{k+1} \subsetneqq \sC_k$.
Furthermore, $I_x\sC_k + \sC_{k+1}$ is an ideal.  (Multiplication of an
element of $I_x\sC_k + \sC_{k+1}$ by a monomial of the form
$f_nU^n$ with $n \geq 1$ clearly yields an element in $\sC_{k+1}$.
Multiplication  by an element of $C_0(X)$
(viewed as an element of $\sP$) produces a product again in
$I_x\sC_k + \sC_{k+1}$.  It follows that $I_x\sC_k + \sC_{k+1}$
satisfies the ideal property.)

We can easily identify $I_x\sC_k + \sC_{k+1}$ as the closed ideal
generated by polynomials of the form
$f_kU^k+\dots+f_nU^n$ where $f_k(x)=0$.

If $\sK$ is any ideal with
$I_x\sC_k + \sC_{k+1} \subset \sK \subseteq \sC_k$, then all elements
of $\sK$ can be written in the form $gU^k+p$, where $p \in \sC_{k+1}$.
For some element of $\sK$, $g(x) \neq 0$.  It now follows that the
coefficients $g$ of $U^k$ in elements of $\sK$ separate points, whence
$\sK = \sC_k$.  Thus, each ideal $I_x\sC_k + \sC_{k+1}$ is maximal
between $\sC_{k+1}$ and $\sC_k$.

On the other hand, let $\sK$ be a maximal ideal between $\sC_{k+1}$
and $\sC_k$.  Again, any element of $\sK$ can be written in the form
$gU^k+p$ with $p \in \sC_{k+1}$ and there must be some element
$x \in X_k$ such that $g(x) = 0$ for all such $g$. (Otherwise, the
Stone Weierstrass theorem again implies that $\sK = \sC_k$.)  So
$\sK \subseteq I_x\sC_k + \sC_{k+1}$.

Thus, $\{I_x\sC_k + \sC_{k+1} \mid x \in X_k \}$ is exactly the family
of maximal ideals which lie between $\sC_{k+1}$ and $\sC_k$.
Essentially the same argument also shows that
$\{\sC_k J_y + \sC_{k+1} \}$ is also the family of maximal ideals
which lie between $\sC_{k+1}$ and $\sC_k$.  In particular, for each
$x \in X_k$ there is one, and only one, $y \in X_{-k}$ such that
$I_x\sC_k + \sC_{k+1} = \sC_k J_y + \sC_{k+1}$.  We next show that
for such a pair $x$ and $y$,  $x = \alpha_k(y)$.

Indeed, let $f_k \in C_0(X_k)$ and $f \in J_y$.  Then
$f_kU^kf = [(f_k \circ \alpha_k)f]\circ \alpha_{-k}U^k$.
For $z \in X_k$,
$[(f_k \circ \alpha_k)f]\circ \alpha_{-k}(z) = f_k(z)f(\alpha_{-k}(z))$.
This vanishes if $\alpha_{-k}(z) = y$, i.e., if $z = \alpha_k(y)$,
showing that all such terms $f_kU^kf \in I_{\alpha_k(y)}$.
  Consequently,
$I_x\sC_k + \sC_{k+1} = \sC_k J_y + \sC_{k+1}$ if, and only if,
$x = \alpha_k(y)$.

Now let $\alpha$ be a partial action of $\bbZ^+$ acting on $X$ and
$\beta$ a partial action acting on $Y$.  Let
$A = C_0(X) \times_{\alpha} \bbZ^+$ and
$B = C_0(Y) \times_{\beta} \bbZ^+$.  Suppose that
$\psi \colon A \to  B$ is a continuous isomorphism.  Let
$\sC_k$, $k=1,2,\dots$ be the ideals considered above for $A$ and
let $\sD_k$ be the corresponding ideals for $B$.  Since $\psi$ is a
continuous isomorphism, $\psi(\sC_k) = \sD_k$, for all $k$.

Now, $\psi$ induces an isomorphism of $A/\sC_1$ onto
$B/\sD_1$ by $p+\sC_1 \mapsto \psi(p) + \sD_1$.  But
$A/\sC_1 \cong C_0(X)$ and $B/\sD_1 \cong C_0(Y)$,
so $\psi$ induces an isomorphism of $C_0(X)$ onto $C_0(Y)$.  In fact,
we can view $C_0(X)$ as a subalgebra of $A$ and $C_0(Y)$ as a
subalgebra of $B$ in a natural way and the isomorphism of $C_0(X)$
onto $C_0(Y)$ is just $\psi$ restricted to these subalgebras.
Consequently, there is a homeomorphism $\tau \colon X \to Y$ such that
$\psi(f) = f \circ \tau^{-1}$, for all $f \in C_0(X)$.

For each $k$, the isomorphism $\psi$ maps $\sC_k$ onto $\sD_k$.  It
follows that $\psi$ carries the ideals which are maximal between
$\sC_{k+1}$ and $\sC_k$ to the ideals which are maximal between
$\sD_{k+1}$ and $\sD_k$.   Suppose that
$\psi(I_x \sC_k+\sC_{k+1}) = I_z\sD_k + \sD_{k+1}$, where
$x \in X_k$ and $z \in Y_k$.

If $f \in C_0(X_k)$ with $f \in I_x$, $p \in \sC_k$, and
$q \in \sC_{k+1}$, then
\[
\psi(fp+q) = \psi(f)\psi(p)+\psi(q) =
f\circ \tau^{-1}\psi(p)+\psi(q)
\]
with $\psi(p) \in \sD_k$ and $\psi(q) \in \sD_{k+1}$.
It follows that $f(\tau^{-1}(z)) = 0$ for all $f \in I_x$, so
$\tau^{-1}(z) = x$.  Thus
\[
\psi(I_x \sC_k + \sC_{k+1}) = I_{\tau(x)}\sD_k + \sD_{k+1}.
\]

In particular, $\tau(X_k) = Y_k$.  Similarly,
$\tau(X_{-k}) = Y_{-k}$ and
$\psi(\sC_k J_y + \sC_{k+1}) = \sD_k J_{\tau(y)} + \sD_{k+1}$.
But
$\sC_kJ_y+\sC_{k+1} = I_{\alpha_k(y)}\sC_k + \sC_{k+1}$ and
$\sD_kJ_{\tau(y)}+\sD_{k+1}=I_{\beta_k(\tau(y))}\sD_k+\sD_{k+1}$.
Since $\psi$ maps the first of these ideals onto the second,
we have
$\tau(\alpha_k(y)) = \beta_k(\tau(y))$, for all $y \in X_{-k}$.

Thus, $\tau \circ \alpha_k = \beta_k \circ \tau$ as a partial
homeomorphism; i.e., the diagram
\[
\begin{diagram}
\node{X_{-k}} \arrow{e,t}{\tau} \arrow{s,r}{\alpha_k}
\node{Y_{-k}} \arrow{s,r}{\beta_k}\\
\node{X_k} \arrow{e,t}{\tau} \node{Y_k}
\end{diagram}
\]
commutes.  This shows that $\alpha$ and $\beta$ are conjugate.

The converse is routine.
\end{proof}

\providecommand{\bysame}{\leavevmode\hbox to3em{\hrulefill}\thinspace}
\providecommand{\MR}{\relax\ifhmode\unskip\space\fi MR }
% \MRhref is called by the amsart/book/proc definition of \MR.
\providecommand{\MRhref}[2]{%
  \href{http://www.ams.org/mathscinet-getitem?mr=#1}{#2}
}
\providecommand{\href}[2]{#2}

\end{document}